\documentclass{amsart}
\usepackage{latexsym}

\usepackage{amssymb,amscd}
\usepackage{amsmath}
\usepackage{graphics}

\usepackage[usenames]{color}

\newtheorem{thm}{Theorem}
\newcommand{\bthm}{\begin{thm}}
\newcommand{\ethm}{\end{thm}}
\newtheorem{cj}{Conjecture}
\newcommand{\bcj}{\begin{cj}}
\newcommand{\ecj}{\end{cj}}
\newtheorem{cor}{Corollary}
\newcommand{\bcor}{\begin{cor}}
\newcommand{\ecor}{\end{cor}}
\newtheorem{prop}{Proposition}
\newcommand{\bprop}{\begin{prop}}
\newcommand{\eprop}{\end{prop}}
\newcommand{\bpf}{\noindent{\bf Proof}}
\newtheorem{defn}{Definition}
\newcommand{\bdf}{\begin{defn}}
\newcommand{\edf}{\end{defn}}
\newtheorem{lem}{Lemma}
\newcommand{\blm}{\begin{lem}}
\newcommand{\elm}{\end{lem}}
\newtheorem{rmk}{Remark}
\newcommand{\brk}{\begin{rmk}}
\newcommand{\erk}{\end{rmk}}

\newtheorem{ex}{Example}
\newcommand{\bex}{\begin{ex}}
\newcommand{\eex}{\end{ex}}
\newcommand{\pt}{\partial}

\newcommand{\T}{{\mathcal T}}

\newcommand{\R}{{\mathbb R}}
\newcommand{\sm}{\setminus}

\newcommand{\diver}{{\rm div}}
\newcommand{\rot}{{\rm rot}}
\newcommand{\sgn}{{\rm sgn}}
\newcommand{\lk}{{\rm lk}}
\newcommand{\Vol}{{\rm Vol}}
\newcommand{\pd}{\partial}

\begin{document}
\title{Asymptotic linking of volume-preserving actions of ${\mathbb R}^k$}

\author{José Luis Lizarbe Chira}
\author{Paul A. Schweitzer, S.J.}

\address{Departamento de Ci\^encias Exatas (VCE), Polo Universit\'ario de Volta Redonda (PUVR)
\\ Universidade Federal Fluminense
\\ Avenida dos Trabalhadores 420, Vila Santa Cecília, Volta Redonda, Rio de Janeiro 27255-125, Brazil.}

\email{jllchira@id.uff.br, chira67@gmail.com}

\address{Departamento de Matematica\\
         Pontificia Universidade Cat\'olica do Rio de Janeiro\\
         Rua Marquês de S\~ao Vicente 225, Gavea, Rio de Janeiro 22453-900, Brazil. }

\email{paul37sj@gmail.com}

\thanks{The first author thanks the CAPES postdoc program 2015, Brazil for support.}

\begin{abstract} We extend V. Arnold's work on asymptotic linking
for two volume preserving flows on a domain in ${\mathbb R}^3$ and
$S^3$ to volume preserving actions of ${\mathbb R}^k$ and
${\mathbb R}^\ell$ on certain domains in ${\mathbb R}^n$ 
and also to linking of a volume preserving action of
${\mathbb R}^k$ with a closed oriented singular $\ell$-dimensional
submanifold in ${\mathbb R}^n$, where $n=k+\ell+1$.

\end{abstract}
\maketitle

\tableofcontents

\section{Introduction}\label{intro} 

V.I. Arnold, in his paper ``The asymptotic Hopf invariant and its
applications" \cite{A} published in 1986 (also see \cite{AK, K, V, CI}), considered a compact
domain $\Omega$ in ${\mathbb R}^3$ or $S^3$ with a smooth boundary and
trivial homology and two divergence free vector fields $X$ and $Y$
in $\Omega$ tangent to the boundary $\partial\Omega$. He defined
an asymptotic linking invariant ${\rm lk}(X,Y)$ that measures the
average linking of trajectories of $X$ with those of $Y$, and
another invariant $I(X,Y)=\int_\Omega\alpha\wedge d\beta$, where
$d\alpha=i_X{\rm \omega}$ and $d\beta=i_Y{\rm \omega}$ (interior
products with the volume form $\omega$ on $\Omega$), and showed
that ${\rm lk}(X,Y)=I(X,Y)$. We extend these results to volume-preserving actions $\Phi$
and $\Psi$ of ${\mathbb R}^k$ and ${\mathbb R}^\ell$ on a compact convex domain $\Omega$ 
with smooth boundary in ${\mathbb R}^n$, where $\Phi$ and
$\Psi$ are tangent to $\partial \Omega$ and $k+\ell=n-1$. 

Arnol'd defines the invariant ${\rm lk}(X,Y)$ as follows. For $p\in
\Omega$ and $T>0$, let $\vartheta_X(p,T)=\{\phi^X_t(p)|0\leq t\leq
T\}$ be the segment of orbit beginning at $p$ and continuing for a
time $T$, and let $\tilde \vartheta_X(p,T)$ be this curve closed
by adding a short path in $\Omega$ from $\phi^X_T(p)$ to $p$.
Define $\tilde \vartheta_Y(q,S)$ similarly. The asymptotic
linking invariant of $X$ and $Y$ is
$${\rm lk}(X,Y)=\int_{\Omega\times\Omega}\tilde{\rm lk}(p,q)
$$ where
$$\tilde{\rm lk}(p,q)=\lim_{S,T\rightarrow\infty}\frac{1}{ST}
{\rm lk}(\tilde \vartheta_X(p,T),\tilde \vartheta_Y(q,S)).$$ Then ${\rm lk}(X,Y)$
is well-defined, since ${\rm lk}(\tilde \vartheta_X(p,T),\tilde
\vartheta_Y(q,S))$ is defined and the limit exists for almost all
$(p,q)\in \Omega\times\Omega$, and furthermore the function $\tilde{\rm
lk}(p,q)$ is in $L^1(\Omega\times\Omega)$ \cite{V}.

The way that Arnol'd closes up partial orbits with short curves was used earlier on by
Schwartzman to define asymptotic cycles for a continuous  flow $\phi$
on a compact polyhedron $X$ \cite{S1}.
Let $\vartheta_\phi(p,T)$ be the partial orbit from $p\in X$ to $\phi_T(p)$,
let $\tilde\vartheta_\phi(p,T)$ be a (possibly singular) loop formed by adding
a short curve, and let $[\tilde\vartheta_\phi(p,T)]\in H_1(X;\R)$ be its first real homology class.
Then the $p$ asymptotic cycle is the limit
$$A_p = \lim 1/t[\tilde\vartheta_\phi(p,T)] \in H_1(M,{\mathbb R})$$
which exists for almost all points $p\in X$, as described in a geometric interpretation 
(\cite{S1}, p. 275). Schwartzman's proof is
quite different, since he uses homomorphisms from the cohomology to $\R$ to define $A_p$.
If the short curves are chosen in a measurable fashion for a normalized invariant measure $\mu$,
then the $\mu$ asymptotic cycle is defined to be the integral 
$A_\mu = \int_X A_p d\mu \in H_1(X;\R)$,
the average of the cycles $A_p$. 

In \cite{S2}, Schwartzman also defines asymptotic cycles for a smooth action
of $\R^k$ on a compact smooth manifold $M^n$. This asymptotic cycle could also be defined by capping off
the boundary of a partial orbit by a small (possibly singular) manifold, if that can be done in a measurable way,
as in the present paper, though this is not carried out in \cite{S2}.

In \S \ref{defsandresults} we define an asymptotic
linking invariant ${\rm lk}(\Phi,\Psi)$ which measures the degree
of linking between orbits of the actions $\Phi$ and $\Psi$ and another
invariant $I(\Phi,\Psi)$ defined in terms of differential forms.
Our main result, Theorem \ref{mainthm} (proven in \S
\ref{link2actions}), states that ${\lk}(\Phi,\Psi)=I(\Phi,\Psi)$. Analogous results are given for the
asymptotic linking of the action $\Phi$ with a closed oriented
$\ell$-dimensional submanifold $N$ (Theorem \ref{thm4}, proven in
\S \ref{linkactionsubman}). 

We use extensions of the gradient, curl,  and divergence to multivectors
in higher dimensions that are presented in \S \ref{gradetc}, and 
in \S \ref{secBS} an extension to higher dimensions of the classical
Biot-Savart formula that gives an inverse for the curl of a
divergence-free vector field on a compact domain ${\mathbb R^3}$. 
A version of the ergodic theorem due to Tempelman \cite{T} that is
used in the proofs is given in \S \ref{ergthmsect}. 

As an application, we show that our invariant gives a lower bound for
the energy of an action in \S \ref{energy}. Examples in which the
invariant is non-trivial are given in the last section, \S
\ref{examples}.

These results are taken from the doctoral thesis \cite{L} of the
first author, under the direction of the second author at the Pontifical Catholic University of Rio de Janeiro (PUC-Rio).
Some similar results were obtained by Garc\'\i a-Comp\'ean and Santos-Silva   in \cite{GS}.
It would be interesting to extend these results to $S^n$ and other
Riemannian manifolds and also to linking of ${\mathbb
R}^k$-actions with leaves of foliations endowed with an invariant transverse
volume form (see \cite{KV}).

\section{Definitions and statements of results}\label{defsandresults} 

Throughout the paper $M$ is an oriented Riemannian $n$-dimensional manifold and
$\Omega\subset M$ is a compact convex domain with smooth boundary
$\partial \Omega$. In the main results of this paper, $M$ will be ${\mathbb
R}^n$ with the standard metric, but many of the details are valid more generally.
We consider a smooth ($C^\infty$) action  
$$\Phi: {\mathbb R}^k\times
\Omega\rightarrow \Omega,$$ 
of the $k$-dimensional real vector
space ${\mathbb R}^k$ on $\Omega$. Then $\Phi$  is defined by $k$ vector fields tangent
to $\partial \Omega$, $X^1, X^2,\dots, X^k$, whose corresponding
flows $\phi^1, \phi^2, \dots, \phi^k$
commute with each other, so that for $t=(t_1,\dots, t_k)\in
{\mathbb R}^k$ and $x\in \Omega$,
$$\Phi(t,x) = \phi^1(t_1,\phi^2(t_2, \dots, \phi^k(t_k,x),\dots)).
$$ In other words, if we set $\Phi_{t}=\Phi(t,\cdot)$ and
$\phi^i_{t_i}=\phi^i(t_i,\cdot)$ for each $i$, then $\Phi_{t}
= \phi^1_{t_1}\circ\dots\circ \phi^k_{t_k}$. As usual, $\phi^i$ is related to $X^i$ by the identity
$\frac{\partial}{\partial t} \phi^i(t,x)=X^i(\phi^i(t,x))$ and the commutation of $\phi^i$ and $\phi^j$ is equivalent to the vanishing of the Lie bracket $[X^i,X^j]$.

\vskip 0.5cm

\bdf A (smooth) action $\Phi: {\mathbb R}^k\times \Omega\rightarrow
\Omega$ on $\Omega$ is {\bf conservative} if it is
volume-preserving (i.e., for each $t\in{\mathbb R}^k$,
$\Phi_{t}:\Omega\to\Omega$ preserves the Riemannian volume form on $M$)
and the generating vector fields $X^i$ are tangent to the boundary $\partial\Omega$.
\label{conservative}\edf

 Let $\Phi: {\mathbb R}^k\times \Omega\rightarrow \Omega$ and
$\Psi: {\mathbb R}^\ell\times \Omega\rightarrow \Omega$ be conservative actions on $\Omega$,
$k+\ell+1=n$.
Let $X= X^1\wedge\dots\wedge X^k$ and $Y= Y^1\wedge\dots\wedge
Y^\ell$ be the exterior products of the $k$ vector fields that
generate the action $\Phi$ and the $\ell$ vector fields that
generate $\Psi$, and let $\omega$ be the volume form on $\Omega$.
Denote the differential forms of degree $r$ on $\Omega$ (resp., the forms that
vanish on $\partial\Omega$) by
$E^r(\Omega)$ (resp., $E^r(\Omega,\partial\Omega)$). 
Since $\Omega$ is convex, their deRham cohomology groups $H^{*}(\Omega;{\mathbb R})$ and
$H^{*}(\Omega,\partial \Omega;{\mathbb R})$ vanish for $0<r<n$.
The differential forms
$i_X\omega\in E^{\ell+1}(\Omega,\partial \Omega)$ and $i_Y \omega\in E^{k+1}(\Omega,\partial \Omega)$ given by the interior
products with $X$ and $Y$ vanish on the boundary $\partial \Omega$ 
since $X$ and $Y$ are tangent to the boundary, and  these forms are closed since the actions are
volume-preserving.  Since $\Omega$ is convex, 
they are exact, so there exist differential forms $\alpha\in E^\ell(\Omega,\partial \Omega)$
and $\beta\in E^k(\Omega,\partial \Omega)$ of degrees $\ell$ and $k$, respectively, such that
$d\alpha=i_X{\rm \omega}$ and $d\beta=i_Y{\rm \omega}$. Then we
define the invariant
$$I(\Phi,\Psi)=\int_\Omega\alpha\wedge d\beta,$$
which obviously does not depend on the choice of $\beta$.
Since $d(\alpha\wedge\beta)=d\alpha\wedge\beta + (-1)^\ell \alpha\wedge d\beta$ and both $\alpha$ and $\beta$ vanish on $\partial\Omega$, Stokes' theorem gives the following result.

\blm This invariant satisfies
$$I(\Phi,\Psi)=(-1)^{\ell+1}\int_\Omega d\alpha\wedge \beta
=(-1)^{(\ell+1)(k+1)} I(\Psi,\Phi).$$
Hence it
depends only on the actions
$\Phi$ and $\Psi$, and not on the choice of the differential forms
$\alpha$ and $\beta$.\label{Iwelldefd} \elm 

We shall define an
asymptotic linking number ${\rm lk}(\Phi,\Psi)$ that measures the
degree of linking between orbits of $\Phi$ and $\Psi$.
For sets $T\subset{\mathbb R}^k$ and $Y\subset\Omega$  we set
$\Phi(T,Y)=\{\Phi(t,y)\ |\ t\in T, y\in Y\}$.
 Let $\mathcal{T}_k$ be the set of $k$-rectangles
$$T=[0,T_1]\times \cdots\times [0,T_k],\ \  (T_1,\dots, T_k)\in {\mathbb R}_+^k$$
where ${\mathbb R}_+^k$ is the space of
$k$-tuples of non-negative real numbers, and fix a   point  $\widetilde p\in\Omega$. 
Then we  let  $\theta_\Phi(p,T)$ be
the closed oriented singular $k$-manifold in the domain $\Omega$
\begin{eqnarray*}\theta_\Phi(p,T)=\Phi(T,p)\cup \sigma(p,T) 
\end{eqnarray*}
where
\begin{eqnarray}
\sigma(p,T) = \Phi(\partial T,p)*\widetilde p \label{sigmapT}
\end{eqnarray}
is the cone composed of the geodesic segments
joining each point of $\Phi(\partial T,p)$ to $\widetilde
p$. We construct the closed oriented singular $\ell$-manifold
$\theta_\Psi(q,S)=\Psi(S,q)\cup \sigma'(q,S)$ in
like manner, replacing $T$ by
$S=[0,S_1]\times\cdots\times[0,S_\ell]\in \mathcal{T}_\ell$ for some
$(S_1,\dots,S_\ell)\in {\mathbb R}_+^\ell$, $\Phi$ by $\Psi$, and $\widetilde p$ by
another point $\widetilde q\neq \widetilde p$.

For fixed $T$ and $S$, since the sum of the dimensions of
$\theta_\Phi(p,T)$ and $\theta_\Phi(q,S)$ is $n-1$, the following lemma holds.
It will  be proved in \S \ref{s-lemmas}.

\blm Fix $T\in\T_k$ and $S\in\T_\ell$.
Then for almost every pair $(p,q)\in \Omega\times \Omega$  the singular manifolds
$\theta_\Phi(p,T)$ and
$\theta_\Psi(q,S))$ are disjoint and therefore ${\rm lk}(\theta_\Phi(p,T), \theta_\Psi(q,S))$ is defined.
\label{measurezero}\elm

The set $D(\Phi,\Psi) = \{(p,T,q,S)\in \Omega\times \T_{k}\times\Omega\times\T_{\ell}\ |\ 
\theta_\Phi(p,T)\cap\theta_\Psi(q,S)=\emptyset\}$, where the compact sets
$\theta_\Phi(p,T)$ and $\theta_\Psi(q,S))$ are disjoint, is clearly open, and  since  it
has full measure, it  must be dense, so we have:

\bcor $D(\Phi,\Psi)$ is an open dense set in $\Omega\times \T_{k}\times\Omega\times\T_{\ell}.$
\label{c:opendense}
\ecor

It follows from the Lemma that the function 
$${\rm lk_{T,S}}(p,q):=\displaystyle\frac{1}{\lambda_k (T)\lambda_\ell (S)}{\rm lk}(
\theta_\Phi(p,T), \theta_\Psi(q,S))$$
is defined for almost all pairs $(p,q)\in \Omega\times\Omega$,
where $\lambda_k(T)=T_1\cdots T_k$ and $\lambda_\ell(S)=S_1\cdots S_\ell$ are the Lebesgue measures
on $\mathbb{R}^k$ and $\mathbb{R}^\ell$. The following theorem, proved in \S \ref{link2actions}, affirms that this function is in $L_1(\Omega\times\Omega)$
and permits us to define the linking index for the orbits of $\Phi$ and $\Psi$. We write $T,S \to \infty$
to signify that $\rm{min}\{T_1,\dots T_k,S_1,\dots, S_{\ell}\}\to\infty$.

\bthm Suppose that $\Omega$ is a compact convex domain in ${\mathbb R}^n$. Let
$\Phi: {\mathbb R}^k\times \Omega\rightarrow \Omega$ and $\Psi:
{\mathbb R}^\ell\times \Omega\rightarrow \Omega$ be conservative
actions with $k+\ell+1=n$. Then

1. The limit function $\displaystyle\lim_{T,S \to \infty} {\rm lk_{T,S}}$ exists as a function in $L^1(\Omega\times \Omega)$, i.e., there is an integrable function $\ \widetilde {\rm lk}_{\Phi,\Psi}:\Omega\times\Omega\to{\mathbb R}$ defined almost everywhere such that
$$ \lim_{T,S \to \infty}\int_\Omega\int_\Omega |{\rm lk_{T,S}}(p,q)-\widetilde {\rm lk}_{\Phi,\Psi}(p,q)|\ dpdq=0.$$ 
2. The integral $\displaystyle \int_\Omega\int_\Omega\widetilde {\rm lk}_{\Phi,\Psi} (p,q)dpdq$ is independent of the choice of the distinct points $\widetilde p$ and $\widetilde q$.
	\label{linkactions} \ethm
Then the {\bf asymptotic
linking number} of $\Phi$ and $\Psi$ is defined to be
$${\rm lk}(\Phi,\Psi):= \int_\Omega\int_\Omega\widetilde {\rm lk}_{\Phi,\Psi} (p,q)dpdq$$
Our main theorem is the following.
\bthm
Under the hypotheses of Theorem \ref{linkactions}, the asymptotic linking number and
the invariant $I(\Phi,\Psi)$ coincide, i.e.,
$${\rm lk}(\Phi,\Psi)=I(\Phi,\Psi).$$\label{mainthm} \ethm

\noindent{\bf Linking of an action with a submanifold.} There is a similar theory for asymptotic linking between a (smooth) conservative action $\Phi: {\mathbb R}^k\times \Omega\rightarrow
\Omega$ and a closed oriented singular $\ell$-submanifold
$N\subset \Omega$, where as above $\Omega$ is a compact convex
domain in $n$-dimensional Euclidean space and $n=k+\ell +1$. As before, let $\alpha$ be an $\ell$-form
on $\Omega$ satisfying $d\alpha = i_X\ \omega$ where the vector
fields  $X^1, X^2,\dots, X^k$ generate the action $\Phi$, $X=
X^1\wedge\dots\wedge X^k$ and let $\omega$ be the volume form on
$\Omega$. Then we define
\begin{eqnarray}I(\Phi,N)= \int_N\ \alpha.\label{IPhiN}\end{eqnarray}

By analogy to the previous case of two actions, we can also define
an asymptotic linking number between the action $\Phi$ and $N$. As
before, let $\theta_\Phi(p,T)=\Phi(T,p)\cup \sigma(p,T)$ with the apex of the cone at $\widetilde p\in \Omega\sm N$.
The proof of the following Lemma is analogous to the proof of Lemma \ref{measurezero} and will also be given in \S \ref{s-lemmas}.
\blm
Fix $T\in\T_k$ and let $N'$ be a compact oriented singular $\ell$-submanifold
$N'\subset \Omega$, possibly with boundary.
 Then for almost every point $p\in \Omega$,  $\theta_\Phi(p,T)\cap N'=\emptyset$.
\label{measurezeroN}\elm
\noindent Hence when $N'=N$, $\displaystyle\frac{1}{\lambda_k (T)}{\rm lk}(\theta_\Phi(p,T), N)$ is 
defined for almost all $p\in\Omega$. Furthermore, the limit as $T\to\infty$ exists in $L^1 ( \Omega)$, and  the integral is well-defined:

\bthm Let $\Phi: {\mathbb R}^k\times \Omega\rightarrow \Omega$ be a conservative
action on a compact convex domain $\Omega$ in ${\mathbb R}^n$ and let $N\subset \Omega$ be a smooth closed oriented $\ell$-manifold, with $k+\ell+1=n$. Then

1. The limit function $\displaystyle {\rm lk}_{\Phi,N}(p):=
\lim_{T \to \infty}\frac{1}{\lambda_k (T)} {\rm lk}(\theta_\Phi(p,T), N)$
exists as a function in $L^1(\Omega)$, i.e., there is an integrable function $\ \widetilde {\rm lk}_{\Phi,N}:\Omega\to{\mathbb R}$ defined almost everywhere such that
$$ \lim_{T \to \infty}\int_\Omega |\frac{1}{\lambda_k (T)} {\rm lk}(\theta_\Phi(p,T), N)-\widetilde {\rm lk}_{\Phi,N}(p)|\ dp=0.$$

2. The integral $\displaystyle \int_\Omega\widetilde {\rm lk}_{\Phi,N} (p)dp$ is independent of the choice of the  point $\widetilde p$.
\label{linkactionN} \ethm
\noindent Then we define the asymptotic linking number of $\Phi$ and $N$ to be
$${\rm lk}(\Phi,N):= \int_\Omega \widetilde{\rm lk}_{\Phi,N}(p)\ dp.$$

\bthm Under the hypotheses if Theorem \ref{linkactionN}, the asymptotic linking number and
the invariant $I(\Phi,N)$ coincide, i.e.,
$${\rm lk}(\Phi,N)=I(\Phi,N).$$
\label{thm4} \ethm

Theorem \ref{linkactions} follows from Proposition \ref{L1function} in \S \ref{link2actions}, the proof of Theorem
\ref{mainthm} is given at the end of \S \ref{link2actions},
Theorem \ref{linkactionN}  follows from Proposition \ref{linklim} in \S \ref{linkactionsubman}, and the proof of
 Theorem \ref{thm4} is given in \S \ref{linkactionsubman}.

\section{Higher Dimensional Vector Algebra}\label{hdva} 

We recall vector algebra on an oriented Riemannian $n$-dimensional manifold with metric $g$.
Let $E_{x,r}=\wedge_r T_xM$ be the $r$th exterior power of the
tangent space $T_xM$ at $x\in M$,
with exterior multiplication $\wedge :E_{x,r}\times E_{x,s}\to
E_{x,r+s}$. The elements of $E_{x,r}$ are
called $r$-{\bf vectors} or {\bf multivectors}.
Recall that the Hodge
operator $*:E_{x,r}\rightarrow E_{x,n-r}$ is defined for any
positive orthonormal basis $e_1,\dots,e_n$ of $T_xM=E_{x,1}$ by
setting
$$*(e_{i_1}\wedge\dots\wedge e_{i_r})=e_{j_1}\wedge\dots\wedge
e_{j_{n-r}},$$ if $(i_1,\dots,i_r,j_1,\dots,j_{n-r})$ is a positive
permutation of $(1,\dots,n)$, and extending over $E_{x,r}$ by
linearity and antisymmetry. Then 
\begin{eqnarray}*\circ  * =(-1)^{r(n-r)}{\rm id}: E_r\to E_r.\label{**}
\end{eqnarray}

The inner product given by the Riemannian metric $<\ ,\ >$
on $T_x M$ defines an inner product on $E_{x,r}$; for decomposable multivectors
$u=u_1 \wedge \dots \wedge u_r$ and $v_1 \wedge \dots \wedge v_r$,  
$u\cdot v={\rm det}(<u_i , v_j>).$ 
This inner product extends to
an ${\mathbb R}$-bilinear product
\begin{equation}\cdot :E_{x,r}\times E_{x,s}\rightarrow E_{x,s-r},\ \
(u,v)\mapsto u\cdot v=*(u\wedge *v),\label{dotgeneralized}\end{equation} and there is also a
generalization to ${\mathbb R}^n$ of the classical cross product
on ${\mathbb R}^3$
\begin{equation}\times:E_{x,r}\times E_{x,s}\rightarrow E_{x,n-r-s},\ \
(u,v) \mapsto u\times v= *(u\wedge v),\label{crossgeneralized}\end{equation} 
In particular, $u\cdot v = u\times *v$. 

\bprop Let $u\in E_{x,r}, v\in E_{x,s}$, and $w\in E_{x,m}$ be
multivectors.
\begin{enumerate}
\item $u\times(v \times w)=u \cdot(v \wedge w)$.
\item If $\ r+s+m=n$, then
$$(u \times v)\cdot w  =\ast (u \wedge v\wedge w).$$
\end{enumerate} \label{uxvxw}
\eprop

\noindent{\bf Proof.}
1. $u \cdot(v \wedge w) = *(u\wedge *(v\wedge w)) =
*(u\wedge(v\times w)) = u\times (v\times w)$.

2. Note that $ u \times v$ and $w$ are both in $E_{x,m}$, so $(u\times v)\cdot w
\in E_{x,0}=\mathbb{R}$ and $(u\times v) \cdot w=w \cdot (u\times v)$.
Now $w
\cdot (u \times v) = \ast(w \wedge \ast \ast(u\wedge
v))=(-1)^{m(r+s)}*(w \wedge u\wedge v)=
 *(u \wedge v \wedge w).$ \qed
\medskip

It follows from item 2 of the preceding Proposition that if the vectors
$u,v,$ and $w$ are decomposable, say $u=v_1\wedge\dots\wedge v_r,
\ v=v_{s+1}\wedge\dots\wedge v_{r+s},$ and
$w=v_{r+s+1}\wedge\dots\wedge v_n$ with $v_i=\sum_j a_{ij}e_j$
for a positive orthonormal basis $e_1,\dots,e_n$, then
\begin{equation}(u\times v)\cdot w=
\det(a_{ij})\label{det-expansion}.\end{equation}

\bex As usual, a multi-index $I$ is an ordered subset $I=(i_1,\dots,i_k)$ of $\{1,\dots,n\}$ with $i_1<i_2< \dots <i_k$, and we set $e_I=e_{i_1}\wedge \dots\wedge e_{i_k}$ and $|I| = k$.
\begin{enumerate}
\item For (ordered) multi-indices $I$ and $J$ we have 
$$e_I\times e_J = e_K$$
if $I\cap J = \emptyset$, $K= \{1,\dots,n\} \setminus (I\cup J)$, and the ordered
union $I\cup J\cup K$ is a positive permutation of $(1,\dots,n)$;  but $e_I\times e_J= 0$ if $I\cap J\neq\emptyset$. 
\item  In addition, $$e_I\cdot e_J = (-1)^{|K|(n-|J|)}e_K$$ if $I  \subset J$, $K=J\setminus I$, 
and the ordered union $I\cup K$ is a positive permutation of $J$; furthermore, $e_I\cdot e_J$ vanishes if $I \not\subset J$.
\end{enumerate}\label{formulas}
\eex

\bprop
For vectors $u,\ v_1 ,\dots,v_k \in T_xM^n$, we have
$$u\cdot(v_1 \wedge \dots \wedge v_k )=
(-1)^{(k-1)(n-k)}\displaystyle \sum_{i=1}^k (-1)^{i-1} ( u\cdot v_i )\
v_1 \wedge \dots {\widehat v_i} \dots \wedge v_k.$$ \label{udotwedge}
\eprop

\noindent{\bf Proof.}
Using (2) of Example \ref{formulas} with $e_I=e_{j_i}$ and
 $e_J=e_{j_1}\wedge\dots\wedge e_{j_k}$ with $1\leq i\leq k$, we have $|K|=k-1$ and so
\begin{equation}e_{j_i}\cdot(e_{j_1}\wedge\dots\wedge
 e_{j_k}) = (-1)^{(k-1)(n-k)+i-1}  e_{j_1}\wedge\dots \widehat{e_{j_i}}\dots\wedge e_{j_k}.\label{eji}\end{equation}
Note that expanding $u=\sum_{i=1}^n u_ie_i$ and
$v=v_1\wedge\dots\wedge v_k$ with
$v_i= \sum_{j_i=1}^n v_{ij_i}e_{j_i}$ we obtain
\begin{equation*}
u\cdot v =
\displaystyle \sum_{j_1 , \dots , j_k =1}^n
\sum_{i=1}^k  u_{j_i}( v_{1j_1}\dots
v_{kj_k})( e_{j_i}\cdot(e_{j_1}\wedge \dots \wedge e_{j_k})).
\end{equation*}
so the desired formula follows by substituting (\ref{eji}) and reassembling the terms $u$ and $v_1,\dots,v_k$.
\qed

\bex For vectors $u, v, w$ in $\R^n$, by the definition of the product $\times$ and Proposition
\ref{udotwedge}, $u\times(v\times w) = u\cdot (v\wedge w) = (-1)^n[(u\cdot v)w- (u\cdot w)v]$.
In particular, in $\R^3$ we have the well-known formula $u\times (v\times w)= (u\cdot w)v - (u\cdot v)w$.

\eex

\section{Extensions of Gradient, Curl, and Divergence.} \label{gradetc} 
Let $E_k=E_k(M)$ be the space
of smooth $k$-vector fields on a Riemannian manifold $M$, and
let $E^k=E^k(M)$ be the dual space of differential $k$-forms. The
inner product $(U,V)\mapsto U\cdot V$ on $E_k$ determines an
isomorphism

\begin{equation}j:E_k\rightarrow E^k,\ \ j(U)(V)= U\cdot V.\label{defj}\end{equation}

The {\em interior product} $i:E_k\times E^r \to E^{r-k}, (X,\alpha)\mapsto i_X\alpha$, is defined
$i_X\alpha(Y) = \alpha(X\wedge Y)$ for $Y\in E_{r-k}$.

\blm Let $\omega$ be the positive unit volume form on $M$. Then $i_X\omega = j(*X).$ \label{j*X}\elm
\noindent{\bf Proof.}
Consider $X=e_I$ where $I=(i_1,\dots, i_k)$ and $J= (j_1,\dots,j_\ell)$  are ordered
multiindices such that $(i_1,\dots,i_k,j_1,\dots,j_\ell)$ is a positive permutation of $(1,\dots,n)$, and let
$\eta_1,\dots,\eta_n$ be the basis dual to a local positive orthonormal basis
$e_1\dots,e_n$.
Then  $$i_{e_I}\omega= \eta_{j_1}\wedge \dots\wedge \eta_{j_\ell} = j(e_J) =  j(*e_I)$$ since $e_J=*e_I$. The lemma follows since every $X\in E_k$ is a linear combination of the
elements $e_I$. \qed

\medskip

The duality between $E_k$ e $E^k$ will be expressed using the isomorphism $j$.  For example, the gradient
operator $\nabla$, defined $\nabla f=j^{-1}(df)$ for a smooth
function $f$ on $M$, can be extended to a linear operator $\nabla: E_k\to E_{k+1}$, $\nabla X = j^{-1}dj(X)$.

\noindent We can also extend the curl and divergence to operators
$\rot:E_k\rightarrow E_{\ell}$
and $\diver:E_k\rightarrow E_{k-1}$ by setting

\begin{equation}\rot(X)=(-1)^{(k+1)\ell}*(\nabla X)\end{equation} and

\begin{equation}\diver(X)=(-1)^{(k+1)\ell}*\nabla(*X)\end{equation}
where we always set $\ell=n-k-1$. On ${\mathbb R}^3$ these
definitions coincide with the classical definitions of curl and
divergence for vector fields.

For the rest of this section we suppose that $M= {\mathbb R}^n$
with the canonical basis $\{e_1,\dots,e_n\}$ and the dual basis $\{dx_1 , \dots dx_n \}$.
For a $k$-vector field of the form $X=fe_{i_1}\wedge\dots\wedge e_{i_k}$ where $f$ is a smooth function it is easy to check that 
$$j(fe_{i_1}\wedge\dots\wedge e_{i_k}) = fdx_{i_1}\wedge\dots\wedge dx_{i_k},\quad
\nabla X=(\nabla f)\wedge e_{i_1}\wedge\dots e_{i_k},
$$
and 
\begin{equation}\diver(X)=(-1)^{k}\sum_{s=1}^k(-1)^{s}
\frac{\partial f}{\partial x_{i_s}}e_{i_1}\wedge\dots\widehat{e}_{i_s}\dots\wedge e_{i_k}.\label{divfe}\end{equation}

Recall that a vector field $U=\sum_{i=1}^n u_i e_i$ on ${\mathbb R}^n$ acts on a function $f$ by setting
$U(f)=\ <U,\nabla f>\ =\ \sum_{i=1}^n u_i \frac{\partial f}{\partial x_i}$.
The action of $U$ on a vector field
$V= \sum_{i=1}^n v_i e_i$ is defined by setting
$$U(V)= \sum_{i=1}^n U(v_i) e_i=\sum_{i,j=1}^n u_j \frac{\partial v_i}{\partial x_j}e_i$$ so the Lie bracket can be written
$[U,V]=U(V)-V(U)$.
\bprop
Let $V=V^1 \wedge \dots \wedge V^k$ be the exterior product of
vector fields
$V^1,\dots, V^k$ on ${\mathbb R}^n$. Then
\begin{eqnarray}
\diver(V)=(-1)^k \displaystyle \sum_{i=1}^k (-1)^{i}  \diver(V^i)\ V^1 \wedge \cdots \widehat{V}^{i}\cdots \wedge V^k\nonumber 
\\
 +(-1)^k \sum_{1\leq i<j \leq k} (-1)^{i+j}[V^i, V^j] \wedge V^1 \wedge \dots\wedge \widehat{V}^i\wedge \dots\wedge \widehat{V}^j\wedge \cdots \wedge V^k\nonumber
\end{eqnarray} where $[V^i, V^j]$ is the Lie bracket. \label{divV}
\eprop

\noindent{\bf Proof.}  Note that this is a dual version of the well-known formula for the exterior derivative of a product of
$1$-forms evaluated on vector fields. Let
\begin{equation} V^i=\sum_{\ell=1}^n v^i_{\ell} \ e_{\ell}\label{expandVi}\end{equation} for every $i$, so expanding $V$
we have
\begin{equation}V= \sum_{\ell_1,\dots,\ell_k=1}^n v^1_{{\ell}_1}\dots
v^k_{{\ell}_k} e_{{\ell}_1}\wedge \dots \wedge e_{{\ell}_k}\label{expandv}.\end{equation} Then by (\ref{divfe})

\begin{eqnarray}\diver(V) & = &
\displaystyle (-1)^k\sum_{i=1}^k\sum_{{\ell}_1 , \dots , {\ell}_k =1}^n (-1)^i
  \frac{\partial( v^1_{{\ell}_1}\dots
v^k_{{\ell}_k})}{\partial x_{{\ell}_i}} 
e_{{\ell}_1}\wedge \dots \wedge \widehat e_{{\ell}_i}\wedge \dots e_{{\ell}_k} \nonumber \\
& = &   \displaystyle(-1)^k\sum_{i,j=1}^k\sum_{{\ell}_1 , \dots , {\ell}_k =1}^n(-1)^i
\frac{\partial v^j_{{\ell}_j}}{\partial x_{{\ell}_i}}  v^1_{{\ell}_1}\dots \widehat{v^j_{{\ell}_j}}\dots
v^k_{{\ell}_k}e_{{\ell}_1}\wedge \dots \wedge \widehat e_{{\ell}_i}\wedge \dots e_{{\ell}_k}.\nonumber
\end{eqnarray}
 Since $\diver(V^i)= \sum_{\ell_i=1}^n \frac{\partial v^i_j}{\partial x_{\ell_i}} $,
the terms with $i=j$ give
\begin{eqnarray*}
(-1)^k \displaystyle \sum_{i=1}^k (-1)^{i}  \diver(V^i)\ V^1 \wedge \cdots \widehat{V}^{i}\cdots \wedge V^k
\end{eqnarray*}
while the remaining terms give the second sum in the proposition; in fact,
if $I_{ab}$ with $a<b$ is the sum of the terms with $(i,j)= (a,b)$ and $(i,j)=(b,a)$, then

\begin{align*} 
I_{ab} = &\ \  (-1)^k\sum_{\ell_a,\ell_b=1}^n
(-1)^{a+b}\left(v^a_{\ell_a}\frac{\partial v^b_{\ell_b}}{\partial x_{\ell_a}}e_{\ell_b}
- v^b_{\ell_b}\frac{\partial v^a_{\ell_a}}{\partial x_{\ell_b}}e_{\ell_a}\right)\wedge\\
& \wedge (v^1_{\ell_1}\dots\widehat{v^a_{\ell_a}}\dots\widehat{v^b_{\ell_b}}\dots  v^k_{\ell_k})
e_{\ell_1}\wedge\dots\widehat{e_{\ell_a}}\dots\widehat{e_{\ell_b}}\dots\wedge e_{\ell_k}\nonumber\\
= &\ \  (-1)^{k+a+b}[V^a,V^b]\wedge V^1\wedge\dots\widehat{V^a}\dots\widehat{V^b}\dots\wedge V^k\nonumber
\end{align*}

since
\begin{eqnarray*}
\sum_{{\ell}_a, {\ell}_b =1}^n  \left(  v^a_{{\ell}_a} \frac{\partial v^b_{{\ell}_b}}{\partial x_{{\ell}_a}}e_{{\ell}_b}
- v^b_{{\ell}_b} \frac{\partial v^a_{{\ell}_a}}{\partial x_{{\ell}_b}}e_{{\ell}_a} \right)  = [V^a,V^ b].\nonumber
\end{eqnarray*}
\qed

\bex If $U$ and $V$ are vector fields in $\R^n$, then by Proposition \ref{divV} and the definitions of $\rot$ and $\times$,
$$\rot(U\times V) = \diver(U\wedge V) = (\diver(V))U - (\diver(U))V - [U,V].$$\eex

\bprop Let  $\omega$ be the positive unit volume form.
Given a $k$-vector field $U\in E_k(\Omega)$ and a $k$-form $\alpha\in E^k(\Omega)$ with $0\leq k\leq  n$,   we have:
	\begin{eqnarray}\alpha(U)\omega=  \alpha\wedge i_U\omega\label{alphaY},
\end{eqnarray}
\begin{eqnarray} dj(U)   = i_{\rot(U)}\omega.\label{djomega}
\end{eqnarray}\label{p-alphaU}\eprop

\noindent{\bf Proof.}  If $U=e_{i_1}\wedge\dots\wedge e_{i_k}$ with $i_1<\dots<i_k$
and $\alpha = dx_{j_1}\wedge\dots\wedge dx_{j_k}$ with $j_1<\dots<j_k$, then $\alpha(U)\neq 0$ if and only if the sequences $(i_1,\dots,i_k)$ and $(j_1,\dots,j_k)$ coincide, and then $\alpha(U)\omega=\omega=\alpha\wedge i_U\omega$. If the two sequences do not coincide, then both sides vanish. By expanding any $U$ and $\alpha$ and using linearity, we conclude that the equation (\ref{alphaY}) holds in general.

Next, $dj(U) = j(\nabla U) = j((-1)^{(k+1)(n-k)}**\nabla U) = j(*\rot(U))$ which is equal to $ i_{\rot(U)}\omega$
by Lemma \ref{j*X},
thus proving (\ref{djomega}).\qed


\section{The Ergodic Theorem for actions of $\mathbb{R}^k$
\label{ergthmsect}} 

In this section we present Theorem \ref{ergthm}, a special case
of Tempelman's
version of the Ergodic Theorem \cite{T} (also see \cite{T1}), for volume-preserving actions of
${\mathbb R}^k$. This result is an essential step in showing that the asymptotic linking invariant 
is well-defined.

Let $M$ be a compact Riemannian manifold (possibly with boundary) with Riemannian volume
form $\mu$ and let $\Phi: \mathbb{R}^k \times M\rightarrow M$ be a conservative action of $\mathbb{R}^k$ on $M$. Let $L^1(M)$
denote the space of measurable real functions $f: M\rightarrow
{\mathbb R}$ such that $\int_M |f|d\mu < \infty$. Consider a
sequence of $k$-rectangles
$$T_n:=[0,T_n^1] \times\dots \times [0,T_n^k],\ n\in {\mathbb N}$$
with each $T^i_n >0$, such that for each $i\ (1\leq i \leq k)\ {\rm
lim}_{n\rightarrow\infty} T_n^i = \infty$. For a function $f\in
L^1(M)$, define a sequence of means $f_n \in L^1(M),\ n\in
{\mathbb N},$ by setting
\begin{eqnarray*}
f_n(p)&:=& 
\frac{1}{\lambda(T_n)}\int_{t\in
		T_n} f(\Phi_{t}(p))d\lambda(t)	\\ 
&=& \frac{1}{T_n^1T_n^2\dots T_n^k} \int_0^{T_n^k} \int_0^{T_n^{k-1}}\dots
\int_0^{T_n^1} f(\Phi_{(t_1,\dots t_k)}(p))dt_1dt_2\dots dt_k \nonumber
\end{eqnarray*}
where $\lambda$ is the Lebesgue measure on $\mathbb{R}^k$ and
$t = (t_1,\dots t_k)$. The following theorem is a special
case of Theorem 6.2 of Tempelman \cite{T} and also of Theorem 3.3 of Lindenstrauss \cite{Lin}.

\bthm (Ergodic Mean Theorem) There is a unique function $\widetilde f$
in $L^1(M)$ to which the sequence $\{f_n\}_{n\in {\mathbb N}}$
converges almost everywhere, i.e.,
\begin{eqnarray*}{\rm lim}_{n\rightarrow \infty} \int_M |f_n - \widetilde f|\ d\mu =0.\end{eqnarray*}
Furthermore, $\widetilde f$ is independent of the choice of the
sequence $\{T_n\}_{n\in {\mathbb N}}$ and satisfies
	\begin{eqnarray*}\int_M \widetilde f\ d\mu = \int_M f\ d\mu.\end{eqnarray*}
\label{ergthm}
 \ethm
Of course, uniqueness of $\widetilde f$ is understood in the sense of $L^1$, i.e.,
two such functions $\tilde f$ agree outside of a set of measure zero.
\bigskip

Lindenstrauss' Theorem 3.3 implies this theorem since ${\mathbb R}^k$ is an amenable
group and $\{T_n\}$ is a tempered F\o lner sequence.

\bigskip

{\noindent \bf Outline of the Proof.} 
First we observe that for a fixed sequence $\{T_n\}$ of $k$-rectangles the
set of $f\in L^1(M)$ for which the Theorem holds is a closed vector subspace of $L^1(M)$.
 Then the essential idea is Tempelman's decomposition of
$L^1(M)$ into invariant functions and functions with zero mean
(Theorem 5.1 of \cite{T}).
Let $W$ be the
vector subspace of $L^1(M)$ generated by functions $h - h\circ\Phi_{ t}$
where $h=\chi_A$ is the characteristic function of a measurable set $A$ and $t\in {\mathbb R}^k$, and let
${\overline W}$ be its closure in $L^1(M)$. One shows that the conclusions of the Theorem hold for $f=h - h\circ \Phi_{t}$, if $h$ 
is the characteristic function of a measurable set $A$ in $\Omega$, 
and consequently for every $f\in W$. By approximation, the same is true for all $f\in \overline W$.

On the other hand, let $I\subset L_1(M)$ be the set of invariant functions where
$f\in L_1(M)$ is invariant if
there exists a measurable set $A$ with
$\mu(M\setminus A)=0$ such that for every $x\in A$ and $t =(t_1,\dots,t_k)
\in {\mathbb R}^k$ we have $f(\Phi_{t}(x)) = f(x)$.
For every invariant function $f$ it is clear that $f_n=f$, so it is easy to see that the conclusions of the
Theorem hold for every $f\in I$ by setting $\widetilde f=f$. Since by Theorem 5.1 of
 \cite{T} every function $f\in L^1(M)$ can be uniquely
represented as a sum
$f = f_1 + f_2$ with $f_1\in I$ and $f_2 \in \overline W$, the Theorem holds for every $f\in  L^1(M)$. \qed

\section{The Generalized Gauss Divergence Theorem for a Multivector Field} \label{generalizedGauss} 

In this section, $\Omega$ is a compact domain with smooth boundary in ${\mathbb  R}^n$.
We define the integral of a $k$-vector field $X= \sum_{1\leq i_1 < \dots <i_k \leq n}\ f_{i_1  \dots i_k } \ e_{i_1}\wedge \dots \wedge e_{i_k}\in E_k(\Omega)$  to be the $k$-vector
\begin{eqnarray}\int_{\Omega} X \ \omega:=\displaystyle \sum_{1\leq i_1 < \dots <i_k \leq n}\ \left(\int_{\Omega}\ f_{i_1  \dots i_k }\ \omega \right) \ e_{i_1}\wedge \dots \wedge e_{i_k}\in E_k(\Omega) \label{defint}\end{eqnarray}
where $\omega$ is the unit volume form.
Using this definition of the integral, we can extend
the Gauss divergence theorem to $k$-vector fields on $\Omega$ with $k>1$.

\bthm (Generalized Gauss Divergence Theorem for a Multivector Field) If $V\in E_k(\Omega)$,
then $$ \int_{\Omega} \diver(V)\ \omega = (-1)^{(k+1)\ell} \int_{\partial \Omega} N\cdot V\ dA\ ,$$
where $N$ is the unit normal vector field pointing outwards along $\partial \Omega$,
$N\cdot V$ is the extended dot product (\ref{dotgeneralized}),
 $\omega$ and $dA$ are the positive unit volume forms on $\Omega$ and ${\partial \Omega}$, and $\ell=n-k-1$. \label{Gaussthm}
\ethm

\noindent{\bf Proof.}  Since every element of $E_k(\Omega)$ is a sum of  decomposable ones, it suffices to prove the proposition for a decomposable $k$-vector $V=V^1 \wedge \dots \wedge V^k$ where $V^i$ is given by (\ref{expandVi}).
Then from (\ref{expandv}) and (\ref{divfe}) we get
$$\diver(V)=(-1)^k\sum_{i=1}^k
\sum_{{\ell}_1 , \dots \widehat{\ell}_i \dots , {\ell}_k =1}^n  (-1)^i
\diver(v^1_{\ell_1}\dots \widehat v^{i}_{\ell_{i}}\dots v^k_{\ell_k}V^i)e_{{\ell}_1}\wedge \dots
\widehat{e}_{{\ell}_i} \dots \wedge e_{{\ell}_k}$$
since $$\diver(v^1_{\ell_1}\dots \widehat v^{i}_{\ell_{i}}\dots v^k_{\ell_k}V^i) =
\sum_{{\ell}_i=1}^n  \frac{\partial( v^1_{{\ell}_1}\dots
v^k_{{\ell}_k})}{\partial x_{{\ell}_i}}.$$
By Stokes' Theorem we have
$$\int_{\Omega} \diver\left( v^1_{\ell_1}\dots \widehat v^{i}_{\ell_{i}}\dots v^k_{\ell_k}\ V^i \right)\omega=\int_{\partial \Omega} v^1_{\ell_1}\dots\widehat v^{i}_{\ell_{i}}\dots v^k_{\ell_k}\ <N, V^i>dA$$
so
\begin{eqnarray*}
\int_{\Omega}\diver(V)\ \omega & = &\displaystyle (-1)^k
\sum_{i=1}^k (-1)^i \int_{\partial \Omega} <N,V^i> V^1 \wedge \dots \widehat{V}^i \dots \wedge V^kdA \\
& = &\displaystyle \int_{\partial \Omega}\left((-1)^k
\sum_{i=1}^k (-1)^i  <N,V^i> V^1 \wedge \dots \widehat{V}^i \dots \wedge V^k \right)dA\\
& = &\displaystyle (-1)^{(k+1)\ell}\int_{\partial \Omega} N \cdot \left( V^1 \wedge \dots \wedge V^k \right)\ dA
\end{eqnarray*}
using $V^j=\sum_{{\ell_j}=1}^n v^j_{\ell_j} e_{\ell_j}$ and Proposition \ref{udotwedge}. \qed

\bcor Set $\Omega - x = \{u-x\in\R^n\ |\ u\in\Omega\}$. For a $k$-vector field $V(x,u)$ on $\R^n\times\R^n$
we have
\begin{eqnarray*}\diver_x\int_{\Omega-x}\hspace{-1em} V(x,u) du = -(-1)^{(k+1)\ell} \int_{\partial\Omega-x}\hspace{-2em}N\cdot V(x,u) dA(u)
+\int_{\Omega-x}\hspace{-1em}\diver_xV(x,u)du. \end{eqnarray*}\label{cor1}
\ecor

\bpf. \ By the change of variables $v=u+x$
$$\int_{\Omega-x} V(x,u) du = \int_\Omega V(x,v-x)dv,$$ 
so

\begin{eqnarray*}\diver_x\int_{\Omega-x}\hspace{-1em} V(x,u) du =
\diver_{1,x}\int_{\Omega} V(x,v-x) dv +
\diver_{2,x}\int_{\Omega} V(x,v-x) dv,
\end{eqnarray*} where the notation indicates that the divergence is calculated with respect to the first or second occurrence of the variable $x$.
Now \begin{eqnarray*} \diver_{1,x}\int_{\Omega} V(x,v-x) dv = &\int_{\Omega} \diver_{1,x}V(x,v-x) dv\nonumber\\
= & \int_{\Omega-x} \diver_xV(x,u) du
\end{eqnarray*}
by reversing the change of variables.
On the  other hand, if we introduce a new variable $z=x$ to separate the two arguments of $V$,
\begin{eqnarray*}
\diver_{2,x}\int_{\Omega} V(x,v-x) dv 
=&\hspace{-6.5em} \diver_x\int_{\Omega} V(z,v-x) dv\\
=&\hspace{-6.5em} -  \diver_v\int_{\Omega} V(z,v-x) dv\\
=&  -(-1)^{(k+1)\ell} \int_{\partial\Omega}N\cdot V(z,v-x) dA(v)\\
=&  -(-1)^{(k+1)\ell} \int_{\partial\Omega-x}N\cdot V(x,u) dA(u)
\end{eqnarray*} by Theorem \ref{Gaussthm}, reversing the change of variables.
Adding the last two expressions gives the desired result. \qed

\section{Extension of the Biot-Savart Formula}\label{secBS}  

We now give an extension of the Biot-Savart formula to higher dimensions.
For a smooth divergence-free vector field $V$ that is tangent to the boundary on a bounded domain $\Omega$ in ${\mathbb  R}^3$, 
it is well known that
the Biot-Savart formula $$BS(V)(x) =\frac{-1}{4\pi}\int\frac{(x-y)\times V(y)}{||x-y||^3}dy$$
gives a right inverse for the curl, i.e., $\rot(BS(V))=V$ (e.g., see \cite{CDTG} \S 5). We generalize this result to $\R^n$
since it will be used in our proofs.

Let $\Omega$ be a bounded domain with smooth boundary $\partial\Omega$ in ${\mathbb R}^n$ and consider $k$
commuting vector fields $V_1,\dots,V_k$ on
$\Omega$ that are divergence-free and tangent to $\partial\Omega$, $1\leq k<n$, with
$\ell=n-k-1$. They generate an action of ${\mathbb R}^k$ on
$\Omega$. Let $V=V_1\wedge\dots\wedge V_k$ be the exterior product of the vector
fields $V_i$.

\bthm For $x\in\Omega$, the $\ell$-vector field
\begin{equation} BS(V)(x)=\frac{(-1)^k}{a_n}\int_\Omega\frac{(x-y)}{||x-y||^n}\times
V(y)dy,\label{eqBS}
\end{equation} where $a_n$ is the $(n-1)$-volume of the unit
sphere in ${\mathbb R}^n$ and we use  the standard Lebesgue measure $dy$ on $\R^n$, satisfies
$${\rm rot}(BS(V))(x)=V(x).$$ \ethm \label{thmBS}

{\noindent \bf Proof.} 
Note that the integral is well defined since the pole along the singular set has order $n-1$.
We  prove the theorem  for $x\in \mathring{\Omega}$, to avoid the problem of a singularity of order $n-1$ when
we integrate along $\partial\Omega$. It will follow by continuity that the theorem holds for every $x\in\Omega$.

By the change of variables $u = y-x$ on  $\Omega-x$, we have
$$BS(V)(x)=\displaystyle\frac{(-1)^{k+1}}{a_n}\int_{\Omega-x} \frac{u}{||u||^n}\times
V(u+x)du.\nonumber$$
Since $\rot(u\times v) = \diver(u\wedge  v)$, from Corollary \ref{cor1}  we get  
\begin{align*}  
I:=& \rot(BS(V))(x)=\frac{(-1)^{k+1}}{a_n} \diver_x\int_{\Omega-x}
\frac{u}{||u||^n}\wedge V(u+x)du = I_1+I_2
\end{align*}
where \begin{eqnarray*}I_1=\displaystyle\frac{(-1)^{k+1}}{a_n} \int_{\Omega-x}\diver_x
\left(\frac{u}{||u||^n}\wedge V(u+x)\right)du\end{eqnarray*}
and
\begin{eqnarray*}I_2= -\displaystyle\frac{(-1)^{k+1+k(\ell+1)}}{a_n}\int_{\partial\Omega-x} 
N\cdot (\frac{u}{||u||^n}\wedge V(u+x))dA(u). \label{I2}\end{eqnarray*}

Applying Proposition \ref{divV} and the facts that
${\rm div}(V^i)=0,\ [V^i , V^j]=0$, and
$\displaystyle \frac{u}{||u||^n}$ does not depend on the variable $x$, we have
$$ I_1=\displaystyle \frac{1}{a_n} \int_{ \Omega-x} \sum_{i=1}^k (-1)^i \left(\frac{u}{||u||^n}\right)_x\left( V^i(u+x)\right)\ V^1 \wedge \dots\widehat{V}^i\dots \wedge V^k (u+x)\ du,$$
where
$\left(\frac{u}{||u||^n}\right)_x\left(V^i \right)$ is the action of the vector field $\frac{u}{||u||^n}$  on $V_i(u+x)$ with derivatives in the variable $x$. 
Expanding the $V_i$'s by (\ref{expandVi}), using the definition of the integral (\ref{defint}), and avoiding 
the singularity at $u=0$, we can write
$$I_1=  \frac{-1}{a_n} \lim_{\epsilon \to 0} \sum_{j_1 , \dots , j_k=1}^{n} \int_{ \Omega'}
\left<\frac{u}{||u||^n}, \nabla_x (v^1_{j_1}\dots v^k_{j_k})(u+x)\right>\hspace{-0.5em}\ du\ e_{j_1}\wedge \dots \wedge e_{j_k} $$
where $\Omega' = (\Omega -x)\setminus \{|| u ||\leq \epsilon \}$ and
$e_{j_i}\wedge e_{j_1}\dots  \widehat{e_{j_i}}\dots e_{j_k}= (-1)^{i-1}e_{j_1}\dots e_{j_k}$.
Now
$$\nabla_x (v^1_{j_1}\dots v^k_{j_k})(u+x)=\nabla_u (v^1_{j_1}\dots v^k_{j_k})(u+x),$$ so,
for $\epsilon>0$ so small that $\{|| u ||\leq \epsilon \}\subset\mathring{\Omega}$, the  integral
$$I(\epsilon):= \int_{ \Omega'}
\left<\frac{u}{||u||^n}, \nabla_x (v^1_{j_1}\dots v^k_{j_k})(u+x)\right>\ du$$
can be written as

\begin{align}
I(\epsilon)&= \int_{ \Omega'}
\left<\frac{u}{||u||^n}, \nabla_u (v^1_{j_1}\dots v^k_{j_k})(u+x)\right>\ du\nonumber\\
&= \int_{ \Omega'}
\left( {\rm div}_u (v^1_{j_1}\dots v^k_{j_k}\frac{u}{||u||^n})-v^1_{j_1}\dots v^k_{j_k}{\rm div}_u (\frac{u}{||u||^n})\right)du\nonumber\\
&= \int_{ \Omega'}
 {\rm div}_u (v^1_{j_1}\dots v^k_{j_k}\frac{u}{||u||^n})du\nonumber\\
 &= \int_{\partial \Omega-x}
 <N, \frac{u}{||u||^n}>\ v^1_{j_1}\dots v^k_{j_k}dA(u)\nonumber\\
&\qquad - \int_{\{|| u ||= \epsilon \}}
  \frac{1}{\epsilon^{n-1}}\ v^1_{j_1}\dots v^k_{j_k}(u+x)\ dA(u)\nonumber
\end{align}  by  Theorem \ref{Gaussthm}, since
$\diver_u(\frac{u}{||u||^n})=0$ on ${\mathbb R}^n$.
Thus $$\lim_{\epsilon\to 0} I({\epsilon})=\int_{\partial \Omega-x}
 <N, \frac{u}{||u||^n}>\ v^1_{j_1}\dots v^k_{j_k}dA(u) - a_n v^1_{j_1}\dots v^k_{j_k}(x),$$ so
 \begin{equation}
 I_1=\displaystyle - \frac{1}{a_n}\int_{\partial \Omega-x}
 <N, \frac{u}{||u||^n}>\ V^1\wedge \dots \wedge V^k (u+x)\ dA(u) + V^1 \wedge \dots \wedge V^k (x).\label{expansionI}
 \end{equation}

Next, returning to (\ref{I2}), we get
\begin{eqnarray*}I_2=& \displaystyle\frac{(-1)^{k\ell}}{a_n}\int_{\partial\Omega-x} 
N\cdot (\frac{u}{||u||^n}\wedge V(u+x))dA(u)\\
=& \displaystyle\frac{1}{a_n} \int_{\partial\Omega-x}< N, \frac{u}{||u||^n}> V(u+x))dA(u). \end{eqnarray*}
by Proposition \ref{udotwedge}, since the $V_i$'s are tangent to $\partial\Omega-x$. Adding the last result to (\ref{expansionI}) we obtain the desired conclusion, $I={\rm rot}(BS(V))(x)=I_1+I_2=V(x).$\qed

\bcor Let $\Omega$ be convex with unit volume form $\omega$ and let
$V\in E_k(\Omega)$ be as above. Then
\begin{eqnarray} dj(BS(V)) =i_V\omega\label{djBS}  \hspace{5em}{\textrm and}
\end{eqnarray}
\begin{eqnarray} I(\Phi,\Psi) = \int j(BS(X))\wedge d\beta
\label{djBSV} \end{eqnarray}\label{c-djBS}
\end{cor}
\bpf. By (\ref{djomega}) and Theorem  \ref{thmBS},												\begin{eqnarray*} djBS(V) = i_{\rot(BS(V))}\omega = i_V\omega. 
\end{eqnarray*}
proving (\ref{djBS}). By Lemma \ref{Iwelldefd}
\begin{eqnarray*} I(\Phi,\Psi) = \int_\Omega \alpha\wedge d\beta
\end{eqnarray*}
is independent of $\alpha$, provided that $d\alpha= i_X\omega$. Then by (\ref{djBS}) with $V=X$
\begin{eqnarray*} I(\Phi,\Psi) = \int j(BS(X))\wedge d\beta.
\qed\end{eqnarray*}

\section{Linking of submanifolds}\label{linksubmans} 

In order to study the asymptotic linking invariant we recall the linking of singular submanifolds in ${\mathbb
R}^n$. Let $N$ and $N'$ be closed, oriented, possibly singular,
disjoint submanifolds of ${\mathbb R}^n$ of dimensions $k$ and
$\ell$, where we always suppose that $n=k+\ell+1$. Then the linking number ${\rm
lk}(N,N')$ of $N$ and $N'$ can be defined as follows. Let $C$ be a compact
oriented singular $k+1$-dimensional manifold in ${\mathbb R}^n$
with $\partial C=N$. By a small deformation of $C$, if necessary,
we may suppose that $C$ is transverse to $N'$ and only intersects
it in non-singular points of $N'$.
Then the linking number of $N$ and $N'$ is defined to be

$${\rm lk}(N,N'):={\sum}_p\ \varepsilon_p$$
where the sum is taken over all points $p\in C\cap N'$, with
$\varepsilon_p= +1$ if the orientation of $C\times N'$ coincides with
that of ${\mathbb R}^n$ or $-1$ if the orientations
are opposite. It is well known that this linking number is
symmetric, does not depend on the choice of $C$, and can also be
calculated as

$${\rm lk}(N,N') = {\rm deg}(f: N\times N'\rightarrow
S^{n-1})$$ where
$$f(p,q):=\frac{q-p}{\|q-p\|}$$ is the normalized vector pointing from $p\in N$ to $q\in N'$
and ${\rm deg}(f)$ is the degree of the mapping $f$ relative to
the orientations of $N$, $N'$, and $S^{n-1}$.
If $N$ and $N'$ are
disjoint images of smooth maps $g: \bar N\rightarrow {\mathbb R}^n$ and
$g': \bar N'\rightarrow {\mathbb R}^n$, then the linking number
can be calculated by
\begin{eqnarray}{\rm lk}(N,N') = \frac{1}{a_n}\int_{\bar N\times\bar N'}\bar f^*(\sigma)\label{lkformula}
\end{eqnarray}
where $\bar f=f\circ(g\times g')$ and $a_n=\int_{S^{n-1}}\sigma$ is the volume form on $S^{n-1}$.


In order to prove the next proposition, we observe that if $(t_1,,t_2,\dots,t_k)$ are local coordinates in $N$,
then the volume form $d\eta$ on $N$ can be written in these coordinates as
\begin{equation} d\eta = \Big|\Big|\frac{\partial}{\partial t_1}\wedge\dots\wedge\frac{\partial}{\partial t_k}\Big|\Big|\ 
dt_1dt_2\dots dt_k. \label{volform}
\end{equation}
and similarly for the volume form $d\eta'$ on $N'$ with local coordinates $s_1,\dots,s_\ell$.

\bprop If $N$ and $N'$ are disjoint immersed closed oriented submanifolds in ${\mathbb
R}^n$, then the linking number ${\rm lk}(N,N')$ can be calculated
by the formula \begin{equation} {\rm
lk}(N,N')=\frac{(-1)^k}{a_n}\int_{p\in N}\int_{q\in N'}
\frac{\big((q-p)\times U(p)\big)\cdot U'(q)}{||q-p||^n}d\eta(p)
d\eta'(q)\label{lkintNN'}\end{equation} where $U(p)$ is a unit
$k$-vector on $N$ at $p$ and $U'(q)$ is a unit $\ell$-vector on
$N'$ at $q$ and $\eta$ and $\eta'$ are the volume measures in $N$
and $N'$.

Furthermore, this formula holds if $N$ and $N'$ are the disjoint images of
smooth manifolds $\bar N$ and $\bar N'$ under smooth singular maps
$g: \bar N\rightarrow {\mathbb R}^n$ and $g': \bar N'\rightarrow
{\mathbb R}^n$, since the images of the singular sets (where
$U(p)=0$ or $U'(q)=0$) have measure zero on $N$ and $N'$, by
Sard's Theorem. \label{linkformula} \eprop

\bpf.\ Note that the  volume form $\sigma  = \sum_{i=1}^n(-1)^{i-1}x_i dx_1\dots\widehat{dx_i}\dots dx_n$
on $S^{n-1}$ can be written $\sigma= i_Ydx_1\dots dx_n$, where $Y = \sum_{i=1}
^nx_ie_i$ is the position vector in $S^{n-1}$. Then, since $dx_1\dots dx_n(Z) = *Z$ for any $Z\in \Lambda_n(\R^n)$,
\begin{eqnarray}
\sigma(v_2\wedge v_3\wedge\dots\wedge v_n) &=& i_Ydx_1\dots dx_n)(v_2\wedge\cdots\wedge v_n)
\nonumber\\
&=& dx_1\dots dx_n(Y\wedge v_2\wedge\dots\wedge v_n)\nonumber\\
&=& *(Y\wedge v_2\wedge\dots\wedge v_n).\label{sigmav}
\end{eqnarray}
On the other hand, using local coordinates $(t_1,\dots, t_k,s_1,\dots,s_\ell)$ in $N\times N'$,
since $f(p,q) =\frac{q-p}{||q-p||}$ and $\bar f= f\circ (g\times g')$, we have 
\begin{eqnarray*}
\frac{\partial \bar f}{\partial t_i}(p,q) = \frac{-1}{||q-p||}\frac{\partial}{\partial t_i}(p)
+ \left[\frac{1}{||q-p||} \right]_{t_i}(q-p),\\
\frac{\partial \bar f}{\partial s_j}(p,q) = \frac{1}{||q-p||}\frac{\partial}{\partial s_j}(q)
+ \left[\frac{1}{||q-p||} \right]_{s_j}(q-p).
\end{eqnarray*}
Setting $\frac{\partial}{\partial t} = \frac{\partial}{\partial t_1}
\wedge\dots\wedge\frac{\partial}{\partial t_k}$
and $\frac{\partial}{\partial s} = \frac{\partial}{\partial s_1}
\wedge\dots\wedge\frac{\partial}{\partial s_\ell}$ we get
\begin{eqnarray}\frac{\partial \bar f}{\partial t}\wedge\frac{\partial \bar f}{\partial s}
&=&\frac{\partial\bar f}{\partial t_1}\wedge\dots\wedge
\frac{\partial\bar f}{\partial t_k}\wedge
\frac{\partial\bar f}{\partial s_1}\wedge\dots\wedge
\frac{\partial\bar f}{\partial s_\ell}\nonumber \\
&=& \frac{(-1)^k}{||q-p||^{k+\ell}}
\frac{\partial}{\partial t}\wedge\frac{\partial}{\partial s}+ W\wedge(q-p)\label{Wform}
\end{eqnarray}
where $W$ is a $(k+\ell-1)$-vector.
Thus,  at the point $(p,q)$ in $N\times N'$ that corresponds to the point 
$\frac{q-p}{||q-p||}\in S^{n-1}$, using local  coordinates and  $k+\ell+1=n$ we get
\begin{eqnarray*}  
\bar f^*(\sigma)(p,q) 
=\sigma(\frac{\partial \bar f}{\partial t}\wedge\frac{\partial \bar f}{\partial s})
dt_1\dots dt_kds_1\dots ds_\ell \\
\end{eqnarray*}\begin{eqnarray*}
&=&  dx_1\dots dx_n\Big(\frac{q-p}{||q-p||}\wedge 
\frac{\partial \bar f}{\partial t}\wedge\frac{\partial \bar f}{\partial s}\Big)
dt_1\dots dt_kds_1\dots ds_\ell\\
&=& \hspace{-0.5em} \frac{(-1)^k}{||q-p||^n}
dx_1\dots dx_n\Big((q-p)\wedge 
\frac{\partial}{\partial t}\wedge\frac{\partial}{\partial s}\Big)
dt_1\dots dt_kds_1\dots ds_\ell \ \hspace{0.5em}{\textrm by\ (\ref{Wform})}\\
&=& \frac{(-1)^k}{||q-p||^n}*\left((q-p)\wedge
\frac{\partial}{\partial t}\wedge\frac{\partial}{\partial s}\right)
dt_1\dots dt_kds_1\dots ds_\ell\hspace{3.5em}
{\textnormal by\ (\ref{sigmav})}  \\
&=& \frac{(-1)^k}{||q-p||^n}*\left((q-p)
\wedge
||\frac{\partial}{\partial t}||U)\wedge 
(||\frac{\partial}{\partial s}||U'\right)
dt_1\dots dt_kds_1\dots ds_\ell \\
&=&\frac{(-1)^k}{||q-p||^n}* \big((q-p)
\wedge U(p)\wedge U'(q) \big)d\eta(p)d\eta'(q)\hspace{5.5em}
{\textnormal  by\ (\ref{volform})}  \\
&=&\frac{(-1)^k}{||q-p||^n}* \big{(}(-1)^{(k+1)\ell}U'(q)
\wedge (q-p)\wedge U(p)\big{)}d\eta(p)d\eta'(q)\\
&=&\frac{(-1)^k}{||q-p||^n}* \big{(}U'(q)
\wedge**[ (q-p)\wedge U(p)]\big{)}d\eta(p)d\eta'(q) \hspace{5.5em}{\textnormal  by\ (\ref{**})}   \\
&=& \frac{(-1)^k}{||q-p||^n}*\big{(}U'(q)
\wedge*[ (q-p)\times U(p)]\big{)}d\eta(p)d\eta'(q)
\hspace{5.3em}{\textnormal  by\ (\ref{crossgeneralized})}\\
&=& \frac{(-1)^k}{||q-p||^n} U'(q)
\cdot \big((q-p)\times U(p)\big)d\eta(p)d\eta'(q)
\hspace{8em} {\textnormal  by\ (\ref{dotgeneralized})}\\
&=& \frac{(-1)^k}{||q-p||^n} \big((q-p)\times U(p)\big)
\cdot U'(q)d\eta(p)d\eta'(q)
\end{eqnarray*}
 since $U'(q)$ and $(q-p)\times U(p)$ are in the same dimension $\ell$ so the dot product commutes.
Thus by (\ref{lkformula})
\begin{eqnarray*} {\rm lk}(N,N')&=&
 \frac{1}{a_n}\int_{p\in N}\int_{q\in N'}\bar f^*(\sigma)(p,q) \\
&=& \frac{(-1)^k}{a_n}\int_{p\in N}\int_{q\in N'}
\frac{\big((q-p)\times U(p)\big)\cdot U'(q)}{||q-p||^n}d\eta(p)
d\eta'(q).  \qed \end{eqnarray*}

\brk (See, e.g., \cite{CDTG})
In dimension $3$, when $N$ and $N'$ are curves parametrized
by arc\-length by $\alpha:[0,t_0]\rightarrow N$ and
$\alpha':[0,s_0]\rightarrow N'$, the formula (\ref{lkintNN'})
becomes the well-known Gauss linking number formula
$${\rm lk}(N,N')=\frac{-1}{4\pi}\int_0^{t_0}\int_0^{s_0}
\frac{\big((\alpha'(s)-\alpha(t))\times \Dot\alpha(t)\big)\cdot
\Dot\alpha'(s)}{||\alpha'(s)-\alpha(t)||^3}dtds.$$ \erk

A double differential form $L(x,y)$ on ${\mathbb R}^n\times
{\mathbb R}^n$ of bidegree $(k,\ell)$, $k+\ell=n-1$, is called a
{\bf linking form} if whenever $N=g(\bar N)$ and $N'=g'(\bar N')$ are disjoint
images of smooth singular maps $g:
\bar N\rightarrow {\mathbb R}^n$ and $g': \bar N'\rightarrow  
{\mathbb R}^n$, where $\bar N$ and $\bar N'$ are 
 closed oriented manifolds of dimensions $k$ and
$\ell$, then we have
$${\rm lk}(N,N') =\int_N\int_{N'}L.$$

\bcor
\begin{equation}
L = L(x,y)=\frac{(-1)^k}{a_n}
\frac{\Big((y-x)\times U(x)\Big)\cdot U'(y)}{||y-x||^n}d\eta(x)
d\eta'(y)\label{linkingform}\end{equation} 
 is a linking form on ${\mathbb
R}^n\times {\mathbb R}^n$, where $U(x)$ is a unit
$k$-vector on $N$ at $x$, $U'(y)$ is a unit $\ell$-vector on
$N'$ at $y$, and $\eta$ and $\eta'$ are the volume measures in $N$
and $N'$. \label{c:linkingform} \ecor

This is evident from Proposition \ref{linkformula}.

\section{Proofs of Lemmas \ref{measurezero} and \ref{measurezeroN}}\label{s-lemmas} 

As in \S \ref{defsandresults}, consider two volume-preserving actions $\Phi: {\mathbb
R}^k\times \Omega\rightarrow \Omega$ and $\Psi: {\mathbb
R}^\ell\times \Omega\rightarrow \Omega$  on a compact convex domain
$\Omega$ in a Riemannian $n$-manifold $M$ tangent to the (smooth) boundary $\partial\Omega$, $n=k+\ell+1$.
Recall that $\T_k$ is the set of $k$-rectangles
$T=[0,T_1]\times\cdots\times[0,T_k]\subset {\mathbb R}^k$ for
$(T_1,\dots, T_k)\in {\mathbb R}_+^k$. 
Fix points $\tilde p, \tilde q\in\Omega, \tilde p\neq\tilde q$,
and consider the geodesic cones $\sigma(p,T),\ (p,T)
\in \Omega\times \T_k,$ and  $\sigma'(q,S),\ (q,S)
\in \Omega\times \T_\ell,$ with apices $\tilde p$ and $\tilde q$, as defined in (\ref{sigmapT}).
We now prove Lemma \ref{measurezero}.

\bigskip

\bpf {\bf\ of Lemma   \ref{measurezero}}.\ 
We must show that for every $T\in\T_k$ and $S\in\T_\ell$
the set $$X = \{(p,q)\in\Omega\times\Omega \ |\ \theta_\Phi(p,T)\cap\theta_\Phi(q,S)\neq\emptyset\}$$ has measure zero in $\Omega\times\Omega$.
Set $$A_q = \Phi(-T,\Psi(S,q)),\ \  B_q=\Phi(-T,\sigma'(q,S)),$$ $$B'_p=\Psi(-S,\sigma(p,T)),\ \  {\rm and}\ \ 
\ C_p=\{q\in\Omega\ |\ \sigma(p,T)\cap\sigma'(q,S)\neq\emptyset\}.$$
Note that for any set $K\subset \Omega$ and $p\in\Omega$, $p\in \Phi(-T,K)\iff \Phi(T,p)\cap K\neq\emptyset$. 
Consequently $$p\in A_q \iff \Phi(T,p)\cap\Psi(S,q) \neq \emptyset,$$
$$p\in B_q \iff \Phi(T,p)\cap\sigma'(q,S)\neq \emptyset, {\rm  and}$$
$$q\in B'_p \iff \Psi(S,q)\cap\sigma(p,T)\neq\sigma(p,T) \emptyset.$$
Since $\theta_\Phi(p,T)=\Phi(T,p)\cup\sigma(p,T)$ and similarly for $\theta_\Psi(q,S)$,
it follows that
$$X = \bigcup_{q\in \Omega}((A_q\cup B_q)
\times\{q\}) \cup \bigcup_{p\in\Omega}(\{p\}\times (B'_q\cup C_p)).$$
Each of the sets $A_p, B_p$, and $B'_q$ is a singular compact $(n-1)$-dimensional submanifold with open
dense complement in $\Omega$, and therefore has measure zero in $\Omega$. 

Next we shall show that if $p\neq\widetilde q$ the set $C_p$ has measure zero in $\Omega$.  
Let $\widetilde N$ be the cone 
consisting of straight segments beginning at $\widetilde q$, passing through a point of
$\sigma(p,T)$, and ending at  a point of $\partial \Omega$. Let $N$ be the closure of the component of
$\widetilde N\setminus  \sigma(p,T)$ that does not contain the point $\widetilde q$.
Now $\sigma'(q,S)$ meets  $\sigma(p,T)$ if and  only if $\Psi(\partial S,q)$ meets $N$.
Thus $C_p = \Psi(-\partial S,N)$, which
is a compact  singular manifold (the product of the image of the union of the  $2\ell$ faces  of $S$
with $N$)  of dimension  $(\ell-1)+(k+1)= n-1$,  so  it has  measure zero.   

Note that each of the sets $$\cup_q(A_q\times\{q\}),\ \cup_q(B_q\times\{q\}),\ \cup_p(\{p\}\times B'_p),\ {\rm and} \ \cup_p(\{p\}\times C_p)$$ is closed and therefore measurable in $\Omega\times\Omega$. Hence the function $f: \Omega\times \Omega\to \{0,1\}$, 
defined by setting $f(p,q)=1$  if $p\in A_q$ and  $0$ otherwise, is measurable. 
Since $A_q$ has measure zero in $\Omega$ for almost all $q\in\Omega$,
 and therefore $\int_\Omega f(p,q)dp = 0$ for almost all $q$, 
Fubini's theorem shows that $\int_{\Omega}\int_{\Omega}f(p,q)dpdq = 0$, which means that the set
$\cup_q(A_q\times\{q\})$ has measure zero in $\Omega\times\Omega$. Parallel arguments show that
the sets $\cup_q(B_q\times\{q\}), \cup_p(\{p\}\times B'_p)$, and $\cup_p(\{p\}\times C_q)$
also have measure zero, so their union $X$ 
has measure zero in $\Omega\times\Omega$, as claimed.
\qed
 
\medskip
 
\bpf {\bf\ of Lemma   \ref{measurezeroN}}.\ The proof is similar to the last proof.
We must show that for every $T\in\T_k$ the set $Y=\{p\in\Omega\ |\ \theta_\Phi(p,T)\cap N'\neq\emptyset\}$ has
measure zero in $\Omega$. Observe that $Y=A\cup C$ where $A=\Phi(T^{-1},N')$ and
$C=\{p\in \Omega\ |\ \sigma(p,T)\cap N')\neq\emptyset\}$. Let $\widetilde B$ be the cone consisting of segments beginning at $\widetilde p$, passing through a point of $N'$, and ending at a point of $\pt \Omega$. Let $B$ be the closure of the component of $\widetilde B\setminus N'$ that does not contain the point $\widetilde p$. As in the previous proof, we find that $C=\Psi(-\pt T,N')$, and then $A$, $C$, and their union $Y$  have measure zero in $\Omega$. \qed


\section{Asymptotic linking of an action and a
submanifold}\label{linkactionsubman} 

Consider a volume-preserving action $\Phi: {\mathbb R}^k\times
\Omega\rightarrow \Omega$ tangent to the boundary on a compact convex domain $\Omega\subset{\mathbb R}^n$ with smooth
boundary and let $N\subset \Omega$ be a closed
singular $\ell$-dimensional oriented submanifold of
$\Omega$, with $k+\ell=n-1$. As before, $\T_k$ is the set of
$k$-rectangles $T=[0,T_1]\times\cdots\times[0,T_k]\subset {\mathbb
R}^k$ for $(T_1,\dots, T_k)\in {\mathbb R}_+^k$,
$\tilde p\in\Omega\setminus N$ is fixed, and $X=X_1\wedge\dots\wedge X_k$ generates $\Phi$. 
According to Lemma 3, for every $T\in\T_k$ 
the sets $\sigma(p,T)$ defined in (\ref{sigmapT}) are disjoint from $N$ for almost all $p\in\Omega$. 
The invariant  $I(\Phi,N) =  \int_N\alpha$ with $d\alpha=i_X\omega$ was defined in (\ref{IPhiN}).

\blm This invariant satisfies $I(\Phi,N) =   \int_NjBS(X)$ and does not depend on the choice of $\alpha$.\label{PhiNindep}\elm
\bpf. By (\ref{djBS})  $djBS(X) = i_X\omega = d\alpha$ so $d(\alpha- jBS(X)) =0$. Since $\Omega$  is convex, $\alpha-jBS(X)$ is exact and there exists a form $\theta$ such that $d\theta=\alpha-jBS(X)$. Then 
$I(\Phi,N) - \int_N jBS(X) = \int_N\alpha-\int_N jBS(X) = \int_N d\theta = \int_{\pt  N}\theta = 0$  since $\partial N = \emptyset$. Clearly $\int_N jBS(X)$ does not depend on $\alpha$.\qed

\medskip

\bprop The following conditions are satisfied:
\begin{enumerate}
\item The sets $\sigma(p,T)$ vary measurably in the sense that
for every $T\in \T_k$ there is a function $h_T:\Omega\rightarrow
{\mathbb R}$ defined by
$$h_T(p)=\frac{1}{T_1\dots T_k}\int_{x\in \sigma(p,T)}\int_{y\in
N}L(x,y),$$ and $h_T\in L^1(\Omega)$, i.e.,
$\int_\Omega|h_T(p)|d\eta(p)<\infty$.
\item The family of functions $\{h_T\}$ converges to zero in
$L^1(\Omega)$, i.e., $$\lim_{T_1,\dots,
T_k\rightarrow\infty}\int_\Omega |h_T(p)|d\eta(p)=0.$$
\end{enumerate}
\label{systsubmans} \eprop

\bpf. To prove (1), let $Y_T := \{p\in\Omega\ |\ \sigma(p,T)\cap N\neq\emptyset\}$ and
note that $h_T(p) = (T_1\cdots T_k)^{-1}\int_{\sigma(p,T)}\int_{N}L(x,y)$ is defined and varies continuously
on the dense open set $\Omega\setminus Y_T$, where the compact sets $\sigma(p,T)$ and $N$ are disjoint.
Then since $Y_T$ has measure zero,
$h_T$ is measurable in $\Omega$. 

To show that $h_T$ is integrable and that the limit converges to zero, we parametrize $\sigma(p,T)$ by setting
$$T^i = [0,T_1]\times\cdots\widehat{[0,T_i]}\cdots\times[0,T_k] $$
and 
\begin{eqnarray}
\pd_{i\delta} T = [0,T_1]\times\dots\times \{t_{i\delta}\}\times\dots\times[0,T_k]\label{bdyT}
\end{eqnarray}
where $t_{i0}=0$ and $t_{i1}=T_i$ are the extremities of the interval $[0,T_i]$.
Then $\pd T= \cup_{i=1}^k\cup_{\delta=0}^1\pd_{i\delta}T$ and
$$\Phi(\pd T,p) = \cup_{i=1}^k\cup_{\delta=0}^1\Phi(\pd_{i\delta}T,p),$$ so
\begin{eqnarray*}\sigma(p,T) = \cup_{i=1}^k\cup_{\delta=0}^1\sigma_{i\delta}(p,T),
\end{eqnarray*} where
$\sigma_{i\delta}(p,T)$ is the cone with base $\Phi(\pd_{i\delta}T,p)$ and apex $\widetilde p$.
It suffices to prove the proposition using each $\sigma_{i\delta}(p,T)$ in place of their union $\sigma(p,T)$.

Parametrize $\sigma_{i\delta}(p,T)$ by 
\begin{eqnarray}\sigma_p(r,t^i) = (1-r)\Phi(t^{i\delta},p)
+ r \widetilde p,\hspace{1em} (r,t^i)\in [0,1]\times \pd_{i\delta} T,\label{sigmap}
\end{eqnarray}
 where $t^i = (t_1,\dots,t_{i-1},t_{i+1},\dots,t_k)$ and $t^{i\delta} = 
(t_1,\dots,t_{i-1},t_{i\delta},t_{i+1},\dots,t_k)$. Then 
$\frac{\pd\sigma_p}{\pd r}(r,t^i) = \widetilde p -\Phi(t^{i\delta},p)$
and $$\frac{\pd\sigma_p}{\pd t_j}(r,t^i) = (1-r)X^i(\Phi(t^{i\delta},p)),$$ where
$X^i=X_1\wedge\cdots \widehat{X_i}\cdots\wedge X_k$.

Hence, setting $|T| = T_1\cdots T_k$ and
$h^{i\delta}_T(p)=\frac{1}{|T|}\int_{x\in \sigma_{i\delta}(p,T)}\int_{y\in
N}L(x,y)$, where $L(x,y)$ is the linking form (\ref{linkingform}), we have
\begin{eqnarray}
|h^{i\delta}_T(p)| &\leq&  \frac{1}{|T|}  \int_{x\in\sigma_{i\delta}(p,T)}\int_{y\in N}|L(x,y)|\nonumber\\
 &=& \frac{1}{|T|}\int_{r\in[0,1]} \int_{t^i\in T^i} 
\int_{y\in N}\widetilde{L}(r,t^i,y,p)dr dt^i d\eta(y)\label{hiT}
\end{eqnarray}
where \begin{eqnarray}\widetilde{L}(r,t^i,y,p) = \frac{\big|[(y-\sigma_p(r,t^i))\times\frac{\pd\sigma_p}
{\pd r \pd t^i}(r,t^i)]\cdot U(y)\big|}{||y-\sigma_p(r,t^i)||^n}\label{tildeL},
\end{eqnarray}
$$\frac{\pd\sigma_p}{\pd r \pd t^i}(r,t^i) = \frac{\pd\sigma_p}{\pd r}\wedge
\frac{\pd\sigma_p}{\pd t_1}\wedge\cdots\widehat{\frac{\pd\sigma_p}{\pd t_i}}\cdots
\wedge\frac{\pd\sigma_p}{\pd t_k},$$
$U(y)$ is the unit $\ell$-vector in $\wedge_\ell(T_y(N))$, and $d\eta(y)$ is the volume measure on $N$. 

\blm There exists a constant $W_i>0$ such that for all  $t^i\in T^i$ and $y\in N$
$$\int_{r\in[0,1]}\int_{p\in\Omega}\widetilde{L}(r,t^i,y,p)dr d\lambda(p) \leq W_i$$
where $d\lambda(p)$ is the euclidean measure on $\Omega$.
\label{Wi} \elm
This lemma will be proven at the end of this section. We use it now to show that 
$h^{i\delta}_T \in L^1(\Omega)$. In fact, by (\ref{hiT}),
\begin{eqnarray*}
\int_{p\in\Omega}|h^{i\delta}_T(p)| &\leq& \frac{1}{|T|}\int_{p\in\Omega}\Big[\int_{r\in [0,1]}
\int_{t^i\in T^i}\int_{y\in N}\widetilde{L}(r,t^i,y,p)dr dt^id\eta(y)\Big]d\lambda(p)\\
&=& \frac{1}{|T|}\int_{t^i\in T^i}\int_{y\in N}\Big[\int_{r\in [0,1]}\int_{p\in\Omega}
\widetilde{L}(r,t^i,y,p)dr d\lambda(p)\Big]dt^id\eta(y)\\
&\leq& \frac{W_i}{|T|}\Big[\int_{t^i\in T^i}dt^i\Big]\Big[\int_{y\in N}d\eta(y)\Big]\\
&=& \frac{W_i\Vol(N)T_1\cdots\widehat{T_i}\cdots T_k}{|T|} = \frac{W_i\Vol(N)}{T_i}.
\end{eqnarray*}
so $h^{i\delta}_T \in L^1(\Omega)$ and
$\lim_{T\to\infty}\int_{p\in\Omega}|h_T(p)|d\lambda(p) = 0$. \qed

\bigskip

\bpf {\bf \ of Lemma \ref{Wi}.} Using $\sigma_p(r,t^i)$ and its derivatives,
\begin{eqnarray*}\frac{\pd\sigma_p}{\pd r \pd t^i}(r,t^i)   
= \frac{\pd\sigma_p}{\pd r}\wedge
\frac{\pd\sigma_p}{\pd t_1}\wedge\cdots\widehat{\frac{\pd\sigma_p}{\pd t_i}}\cdots
\wedge\frac{\pd\sigma_p}{\pd t_k}
\end{eqnarray*}  \begin{eqnarray*}
= (1-r)^{k-1}[\Phi(t^{i\delta},p) - \widetilde p]\wedge X_1(\Phi(t^{i\delta},p)) \wedge\cdots
\widehat{X_i(\Phi(t^{i\delta},p))}\cdots \wedge
 X_k(\Phi(t_{i\delta},p))\\
= (1-r)^{k-1}[\Phi(t^{i\delta},p) - \widetilde p]\wedge X^i(\Phi(t^{i\delta},p))
\end{eqnarray*}
where $X^i = X_1\wedge\cdots \widehat{X_i}\cdots\wedge X_k$.
Note that $(1-r)^{k-1}\leq 1$, $|\Phi(t^{i\delta},p) - \widetilde p|$ is less than or equal to the
diameter $D$ of $\Omega$, 
there is a constant $B$ such that $||X^i(p)|| \leq B$ for all $p\in\Omega$, and $||U(y)||=1$, so by (\ref{tildeL}) we have

\begin{eqnarray*}\widetilde{L}(r,t^i,y,p)
&\leq& \frac{||\frac{\pd\sigma_p}{\pd r\pd t^i}(r,t^i)||\ ||U(y)||}{||\sigma_p(r,t^i)-y||^{n-1}}\\
&\leq& \frac{(1-r)^{k-1}||\Phi(t^{i\delta},p) - \widetilde p||\ ||X^i(\Phi(t^{i\delta},p))||\ ||U(y)||)}
{||\sigma_p(r,t^i)-y||^{n-1}}\\
&\leq& \frac{DB}{||(1-r)\Phi(t^{i\delta},p) + r\widetilde p - y||^{n-1}}.
\end{eqnarray*}
Thus

$$\int_{r\in[0,1]}\int_{p\in\Omega}\widetilde{L}(r,t^i,y,p)dr d\lambda(p)$$ 
\begin{eqnarray}\leq
\int_{r\in[0,1]}\int_{p\in\Omega}\frac{DB}{||(1-r)\Phi(t^{i\delta},p) + r\widetilde p - y||^{n-1}}
dr d\lambda(p).\label{elly}
\end{eqnarray} 
Now for $\widetilde p\notin N$ there exists $\epsilon>0$ such that
for all $y\in N$ and $r\in [1-\epsilon,1]$
$$||(1-r)\Phi(t^{i\delta},p) + r\widetilde p - y|| \geq d/2,$$
where $d$ is the distance from $\widetilde p$ to $N$. Then

$$\int_{r\in[1-\epsilon,1]}\int_{p\in\Omega}\frac{DB}{||(1-r)\Phi(t^{i\delta},p) + r\widetilde p - y||^{n-1}}
dr d\lambda(p)$$
$$\leq \int_{r\in[1-\epsilon,1]}\int_{p\in\Omega}
\frac{DB}{(d/2)^{n-1}}dr d\lambda(p)= \frac{2^{n-1}DB\epsilon}{d^{n-1}}.$$

On the other hand, for $r\in[0,1-\epsilon]$,  $\Phi(t^{i\delta},\cdot) = \Phi_{t^{i\delta}}$ is a volume-preserving diffeomorphism
of $\Omega$, so we can make the substitution $p'=\Phi(t^{i\delta},p)$ and get
$$ \hspace{-8em}\int_{r\in[0,1-\epsilon]}\int_{p\in\Omega}
\frac{DB}{||(1-r)\Phi(t^{i\delta},p) + r\widetilde p - y||^{n-1}}dr d\lambda(p)$$ 
$$\hspace{10em}= \int_{r\in[0,1-\epsilon]}\int_{p'\in\Omega}
\frac{DB}{||(1-r)p' + r\widetilde{p} - y||^{n-1}}dr d\lambda(p').$$
Now for each $r$ we let $p_r = (1-r)p' + r\widetilde p$.
Then $d\lambda(p_r) = (1-r)^nd\lambda(p')$ and $\Omega$ is replaced by
by $\Omega_r \subset \Omega$ (a contraction moving towards $\widetilde p$), so
$$\int_{r\in[0,1-\epsilon]}\int_{p'\in\Omega}
\frac{DB}{||(1-r)p' + r\widetilde{p} - y||^{n-1}}dr d\lambda(p')$$
$$= \int_{r\in[0,1-\epsilon]}\frac{1}{(1-r)^{n}}\int_{p_r\in\Omega_r}
\frac{DB}{||p_r - y||^{n-1}}dr d\lambda(p_r)$$
$$\leq \int_{r\in[0,1-\epsilon]}\frac{1}{\epsilon^{n}}\int_{p_r\in\Omega}
\frac{DB}{||p_r - y||^{n-1}}dr d\lambda(p_r)$$
$$\leq \int_{r\in[0,1-\epsilon]}\frac{DB\Gamma}{\epsilon^n}dr =\frac{DB\Gamma(1-\epsilon)}{\epsilon^n}
$$
since $1-r\geq \epsilon$ and $\Omega_r\subset\Omega$,
by the following lemma, which
 holds since the singularity at $q$ has order $n-1$, and that is less than the dimension $n$.
\blm There is a constant $\Gamma$ such that the function
$$g(q) = \int_{\Omega\setminus\{q\}} \frac{1}{||p-q||^{n-1}}d\lambda(p)$$ satisfies $|g(q)| \leq \Gamma$ 
 for all $q\in \Omega.$ \label{Gamma}
\elm\qed

Combining the last two results with (\ref{elly}), we get
$$\int_{r\in[0,1]}\int_{p\in\Omega}\widetilde{L}(r,t^i,y,p)dr d\lambda(p)
\leq \frac{2^{n-1}DB\epsilon}{d^{n-1}} + \frac{DB\Gamma(1-\epsilon)}{\epsilon^n} =: W_i.$$ \qed

\medskip

Since $\theta_\Phi(p,T)$ and $N$ are disjoint for almost all $(p,T)\in \Omega\times \T_k$, the
linking number ${\rm lk}(\theta_\Phi(p,T),N)$ is
defined on an open dense set. Then we have

\bprop The limit $$ \widetilde{\rm lk}_{\Phi,N}(p)=\lim_{T_1,\dots,
T_k\rightarrow\infty} \frac{1}{T_1\dots T_k}{\rm
lk}(\theta_\Phi(p,T),N)$$ exists as an integrable
$L^1$-function on $\Omega$ and does not depend on the choice of the point
$\tilde p\in\Omega\setminus N$.
\label {linklim}\eprop 

\bpf.  
\begin{eqnarray}\hspace{3em}{\rm lk}(\theta_\Phi(p,T),N) = \int_{\Phi(T,p)}\int_N L + \int_{\sigma(p,T)}\int_N L.\label{2links}
\end{eqnarray}
By Proposition \ref{systsubmans}, 
\begin{eqnarray}\lim_{T\to\infty}\frac{1}{T_1\cdots   T_k} \int_{\sigma(p,T)}\int_N L = 0.\label{zerolimit}
\end{eqnarray}
Let  $$g(p) = \frac{(-1)^k}{a_n}\int_{y\in N}\frac{(y-p)\times X(p)\cdot  U(y)}{||y-p||^n}d\eta(y)$$
where $U$ is the positive unit $\ell$-form on $N$.
The function $g$ is smooth on $\Omega\setminus N$. Then
\begin{eqnarray*}|g(p)| \leq 
\frac{1}{a_n}\int_{y\in N}\frac{||y-p||\hspace{0.1em}||X(p)||\hspace{0.1em}   ||U(y)||}{||y-p||^n}d\eta(y).
\end{eqnarray*}
Let $K$ be an upper bound for $||X(p)||$, $p\in\Omega$. Since $||U(y)||=1$,
$$|g(p)| \leq   \frac{K}{a_n}\int_{y\in N}\frac{1}{||y-p||^{n-1}}d\eta(y).$$
By Fubini's Theorem
\begin{eqnarray*}\int_{p\in\Omega}|g(p)|d\lambda(p) &\leq& \frac{K}{a_n}\int_{y\in  N}\int_{p\in\Omega}
\frac{1}{||y-p||^{n-1}}d\lambda(p)d\eta(y)\\
&\leq& \frac{K\Gamma}{a_n}\int_Nd\eta  = \frac{K\Gamma Vol(N)}{a_n} \end{eqnarray*}
so  $g\in L^1(\Omega)$. On the other hand, note that
$$\int_{x\in\Phi(T,p)}\int_{y\in N}L(x,y)=$$
\begin{eqnarray*}&=& \int_0^{T_1}\dots\int_0^{T_k}\int_{y\in N}
\frac{(y-\Phi(t,p))\times X(\Phi(t,p))\cdot U(y)}{||y-\Phi(t,p)||^n} d\eta(y)dt\\
&=& \int_0^{T_1}\dots\int_0^{T_k} g(\Phi(t,p))dt.
\end{eqnarray*}
Thus, by (\ref{2links}), (\ref{zerolimit}), and the Ergodic  Theorem,  since $g\in L^1(\Omega)$, the limit
$$\lim_{T\to\infty}\frac{1}{T_1\cdots   T_k}{\rm  lk}(\theta_\Phi(p,T),N)
= \lim_{T\to\infty}\frac{1}{T_1\cdots   T_k}\int_0^{T_1}\dots\int_0^{T_k}g(\Phi(t,p))dt
$$  exists and defines an $L^1$ function $\widetilde{\rm lk}_{\Phi,N}(p)$ on $\Omega$ that satisfies
$$\int_{p\in\Omega}\widetilde{\rm lk}_{\Phi,N}(p)d\lambda(p)  = \int_{p\in\Omega}g(p)d\lambda(p)$$
and does not depend on the choice of $\widetilde p$. \qed 

\medskip

Then we define the
asymptotic linking invariant to be ${\rm
lk}(\Phi,N)=\int_\Omega\widetilde{\rm lk}_{\Phi,N}(p)\ d\eta$ and prove Theorem \ref{thm4}, which states that
${\rm lk}(\Phi,N)=I(\Phi,N)$. 

\medskip

\bpf \ {\bf of Theorem \ref{thm4}.}

\begin{eqnarray*}
{\rm lk}(\Phi,N) &=& \int_{p\in\Omega}\widetilde{\rm lk}_{\Phi,N}(p)d\lambda(p)
=\int_{p\in\Omega}g(p)d\lambda(p)\\
&=& \int_{p\in\Omega}\frac{(-1)^k}{a_n}\int_{y\in N}\frac{(y-p)\times X(p)\cdot U(y)}{||y-p||^n}
d\eta(y)d\lambda(p)\\
&=&  \int_{y\in N}\Big[\frac{(-1)^k}{a_n}\int_{p\in\Omega}\frac{(y-p)\times X(p)}{||y-p||^n}d\lambda(p)\Big]
\cdot U(y)d\eta(y)
\end{eqnarray*}
by Fubini's Theorem, so by (\ref{eqBS}) and the definition of the isomorphism $j$ \begin{eqnarray*}
{\rm lk}(\Phi,N) = \int_N BS(X)\cdot Ud\eta
= \int_N jBS(X)(U)d\eta.
\end{eqnarray*}
Then since $U$ is a unit $\ell$-vector and $d\eta$ is a unit $\ell$-form, Lemma \ref{PhiNindep} shows that

 $${\rm lk}(\Phi,N)=  \int_N jBS(X) = I(\Phi,N). \quad\quad\quad\qed$$


\section{Asymptotic linking of two actions\label{link2actions}} 


In this section, we assume that $M=\R^n$, so $\Omega$ is a compact convex region  with smooth boundary in $\R^n$ and consider
volume-preserving actions $\Phi$
and $\Psi$ of ${\mathbb R}^k$ and ${\mathbb R}^\ell$ that are tangent to the boundary  on $\Omega$, $k+\ell=n-1$, as in \S \ref{defsandresults}. 
 Recall that $D(\Phi,\Psi)\subset   \Omega\times \T_{k}\times\Omega\times\T_{\ell}$ is the dense open set
 of points $(p,T,q,S)$ for which $\theta_\Phi(p,T)$ and $\theta_\Psi(q,S))$ are disjoint.

\bprop 
The following conditions are satisfied:
\begin{enumerate}
\item The functions $(p,T)\mapsto\theta_\Phi(p,T)$ and $(q,S)\mapsto\theta_\Psi(q,S)$ are continuous 
functions on $\Omega\times\R^k$.
Furthermore, the function $\int_{\theta_\Phi(p,T)}\int_{\theta_\Psi(q,S)} L(p,q)$ is continuous on 
$D(\Phi,\Psi)$ and therefore  measurable.

\item The limits 
\begin{eqnarray}\lim_{T,S\rightarrow\infty} \frac{1}{\lambda_k(T)\lambda_\ell(S)}\int_{p\in \Omega}\int_{q\in\Omega}\Big\{\int_{A_p}\int_{B_q} L(p,q)\Big\}dpdq=0,\label{limint}\end{eqnarray} where we set
$(A_p,B_q)$ equal to $(\Phi(T,p),\sigma'(q,S)),
(\sigma(p,T),\Psi(S,q))$, and
($\sigma(p,T),\sigma'(q,S))$, exist, and all three limits are zero.
\end{enumerate}
\label{smallnbds} \eprop 

\bpf. (1)  Since the actions are continuous and line segments depend continuously on their extremities,  it is clear that   
the functions $(p,T)\mapsto\theta_\Phi(p,T)$ and $(q,S)\mapsto\theta_\Psi(q,S)$ are continuous,
and so the function $\int_{\theta_\Phi(p,T)}\int_{\theta_\Psi(q,S)} L(p,q)$ is continuous and measurable
on the dense open set $D(\Phi,\Psi)$.

Proof of (2). 
As before, $T,S\to \infty$ means that $\min(T_1,\dots,T_k,S_1,\dots,S_\ell)\to\infty$. When the compact sets $A_p$ and $B_q$ are disjoint, it is clear that the integral $\int_{A_p}\int_{B_q} L(p,q)$ converges, but it is not evident that the integral in (\ref{limint}) converges, although the integrand is measurable.
First, consider $A_p=\Phi(T,p)$ and $B_q=\sigma'(q,S)$. We decompose  $\sigma'(q,S)=\cup\sigma'_{j\varepsilon}(q,S)$
analogous to the decomposition (\ref{bdyT}) of $\sigma(p,T)$ with the parametrization (\ref{sigmap}).
Let $s_{j0}=0$  and $s_{j1}=S_j$ be
the extremities of the interval $[0,S_j]$. Note that $\Psi(\partial S,q)$   is  the  union of  $2\ell$ sets,    
$\Psi(\partial S,q)=\cup_{j=1}^\ell\cup_{\varepsilon=0}^1\Psi(\partial_{j\varepsilon}S,q)$, $\varepsilon\in \{0,1\}$, where
$$\partial_{j\varepsilon}S =   [0,S_1],\times\dots\times\{s_{j\varepsilon}\}\times\dots\times[0,S_\ell],$$
so the  singular submanifold 
\begin{eqnarray}\sigma'(q,S) = \cup_{j=1}^\ell\cup_{\varepsilon=0}^1\sigma'_{j\varepsilon}(q,S)\label{sigmaj}
\end{eqnarray}
where $\sigma'_{j\varepsilon}(q,S)$ is the cone  joining  $\Psi(\partial_{j\varepsilon}S,q)$ to the  vertex $\widetilde q$. 
We shall prove the Proposition for $B=\sigma'_{j\varepsilon}(q,S)$ instead of $\sigma'(q,S)$; then the same proof works for the other components of $\sigma'(q,S)$.  

Let $S^j=[0,S_1],\times\dots\widehat{[0,S_j]} \dots\times[0,S_\ell]$. To each point $$s^j=   (s_1,\dots,s_{j-1},s_{j+1},\dots,s_\ell) \in S^j$$ we naturally  associate the point $s^{j\varepsilon} = (s_1,\dots,  s_{j-1},s_{j\varepsilon},s_{j+1},\dots, s_\ell)\in  \partial_{j\varepsilon}S$. We use the parametrizations $x_p(t)=\Phi(t,p),   t\in T,$ of $\Phi(T,p)$ and  
$$y_q(u,s^j)=  (1-u)\Psi(s^{j\varepsilon},q)  +   u\widetilde q,  \quad (u,s^j)\in [0,1]\times S^j,\label{yq}$$ of 
$\sigma'_{j\varepsilon}(q,S)$. Note  that  $\frac{\pt x_p}{\pt t_i}=X_i,  \frac{\pt  y_q}{\pt s^j}=  (1-u)Y_j$
and $\frac{\pt y_q}{\pt  u} = \widetilde q - \Psi(s^{j\varepsilon},q)$. Since $\Omega$ is compact,
there is a constant $C$ that  is  a common upper bound for
$||X(p)||= ||X_1\wedge\dots\wedge X_k(p)||, |Y^j(q)||=||Y_1\wedge\dots\widehat{Y_j}\dots\wedge Y_\ell(q)||$  and for $|\widehat q   - \Psi(s^{j\varepsilon},q)|, p,q\in \Omega$. Recall that for multivectors  $||(u\times v)\cdot w||\leq ||u||\hspace{0.1em}||v||\hspace{0.1em}||w||$. Then
\begin{eqnarray*}
\Big|\int_{A_p}\int_{B} L(p,q)\Big|  &=& \Big|\int_{\Phi(T,p)}\int_{\sigma'_{j\varepsilon}(q,S)}L\Big|\\
&\leq& \int_{\Phi(T,p)}\int_{\sigma'_{j\varepsilon}(q,S)}|L|,\end{eqnarray*}
but using (\ref{linkingform}), $d\eta(x) = d\lambda(t)$, and $d\eta'(y) = (\widetilde q - \Psi(s^{j\varepsilon},q))d\lambda(s^j)du,$
\begin{eqnarray*}|L| &\leq& \frac{1}{a_n}
\frac{||x_p(t)-y_q(u,s^j)||\hspace{0.1em}||X(\Phi(t,p)||\hspace{0.1em}||Y^j(y_q(u,s^j))||
}{||x_p(t)-y_q(u,s^j)||^n}
d\eta'(y) d\eta(x) \\
&\leq&
C'\int_0^1\frac{1}{||x_p(t)-y_q(u,s^j)||^{n-1}}d\lambda(s^j)du d\lambda(t),\end{eqnarray*}
where   $C'=C^3/{a_n}$, so

\begin{eqnarray*} \Big|\int_{A_p}\int_{B}L\Big|
\leq C'\int_T\int_0^1\int_{S^j}\frac{1}{||x_p(t)-y_q(u,s^j)||^{n-1}}d\lambda(s^j)du d\lambda(t).
\end{eqnarray*}

Integrating  $|\int_{A_p}\int_{B}L|$ on $\Omega\times\Omega$ we have
\begin{eqnarray*}\hspace{-16em} \int_{p\in\Omega}\int_{q\in\Omega}\Big|\int_{A_p}\int_{B}L\Big|d\lambda(p)d\lambda(q)
\end{eqnarray*}
\begin{eqnarray*} &\leq& C'\int_{p\in\Omega}\int_{q\in\Omega}\int_T\int_0^1\int_{S_j}
\frac{1}{||x_p(t)-y_q(u,s^j)||^{n-1}}
d\lambda(s^j)du d\lambda(t) d\lambda(p)d\lambda(q)\\
&\leq& C'\int_T\int_0^1\int_{S_j}\int_{q\in\Omega}\Big[\int_{p\in\Omega}
\frac{1}{||\Phi(t,p)-y_q(h,s^j)||^{n-1}}d\lambda(p)\Big]d\lambda(q)
d\lambda(s^j)du d\lambda(t) 
\end{eqnarray*}
by Fubini's Theorem, since we shall see that the last integral  converges.

Since the action $\Phi_t$ preserves the volume, if we set
$\Phi(t,p) = p'$, the measure $d\lambda(p')$ coincides with $d\lambda(p)$, and the last integral becomes
\begin{eqnarray}\hspace{2em}C'\int_T\int_0^1\int_{S_j}\int_{q\in\Omega}\Big[\int_{p\in\Omega}
\frac{1}{||p'-y_q(u,s^j)||^{n-1}}d\lambda(p')\Big]d\lambda(q)
d\lambda(s^j)du d\lambda(t). \label{fubini}\end{eqnarray}
Lemma \ref{Gamma} shows that this integral coverges. Then,
working backwards, it follows that all the previous 
integrals in this proof also converge. The
integral (\ref{fubini}) is less than or equal to
\begin{eqnarray*}  C'\int_T\int_0^1\int_{S^j}\int_{q\in\Omega}\Gamma  d\lambda(q)d\lambda(s^j)du d\lambda(t)
&\leq& C'\Gamma Vol(\Omega)Vol(T)Vol([0,1])Vol(S^j)\\
&=& C'\Gamma Vol(\Omega)T_1\dots T_kS_1\dots \widehat{S_j}\dots S_\ell.
\end{eqnarray*}
In the limit we have
\begin{eqnarray*}
0 &\leq& \lim_{T,S\to \infty } \frac{1}{T_1\dots T_kS_1\dots S_\ell}\int_{p\in \Omega}\int_{q\in \Omega}
\Big|\int_{A_p}\int_{B}L\Big||d\lambda(p)d\lambda(q)\\
&\leq& \lim_{T_1,\dots,T_k,S_1,\dots,S_\ell\to \infty} \frac{C'\Gamma Vol(\Omega)}{S_j} = 0,
\end{eqnarray*}
so (\ref{limint}) holds for $A_p=\Phi(T,p)$ and $B_q=B =\sigma'_{j\varepsilon}(q,S)$.
Thus the limit vanishes for $\Phi(T,p)$ and $\sigma'(q,S)$ and similarly for the case
$A_p=\sigma(p,T)$ and $B_q=\Psi(S,q)$. 

\bigskip


For the case when $A_p=\sigma(p,T)$ and $B_q=\sigma'(q,S)$, we use
 the decompositions (\ref{bdyT}) of $\sigma(p,T)$ and (\ref{sigmaj}) of  $\sigma'(q,S)$ and the parametrizations 
$$x_p= \sigma_p(r,t^i) =  (1-r)\Phi(t^{i\delta},p)
+ r \widetilde p,\hspace{1em} (r,t^i)\in [0,1]\times T,$$
and
$$y_q = \sigma'_q(u,s^j) =  (1-u)\Phi(t^{j\varepsilon},q)
+ u \widetilde q,\hspace{1em} (u,s^j)\in [0,1]\times T,$$
of $\sigma_{i\delta}(p,T)$ and  $\sigma'_{j\varepsilon}(q,S)$,
with $t^i, t^{i\delta}, s^j$ and $s^{j\varepsilon}$ as before.  
Then we have

\begin{eqnarray*} 
|L(x_p,y_q)| &\leq& \frac{1}{a_n}
\frac{||x_p - y_q||\hspace{0.1em}||X(x_p)||\hspace{0.1em}||Y^j(y_q)||
}{||x_p-y_q||^n}d\eta'(y_q) d\eta(x_p) \\
&\leq&
\frac{C}{||x_p - y_q||^{n-1}}d\lambda(s^j)du d\lambda(t),\end{eqnarray*}
where   $Ca_n$ is an upper bound for $|X(p)|\ |Y(q)|$.

It suffices to show that the limit of
$$L = \frac{1}{|S|\ |T|}
\int_{p\in\Omega}\int_{q\in\Omega}\int_{x_p\in\sigma_{i\delta}(p,T)}\int_{y_q\in\sigma'_{j\varepsilon}(q,S)}
\frac{C}{||x_p - y_q||^{n-1}}d\lambda(s^j)du d\lambda(t)dpdq$$
converges to zero as $S,T\to\infty$. We shall do this in three cases.

\medskip

Case 1. $r,u\in[1-\epsilon,1]$, where $\epsilon >0$ is such that $||y_q-x_p||\geq d/2$ when
$u,r\in[1-\epsilon,1]$ and $d$ is the distance from
$\widetilde p$ to $\widetilde q$. Such an $\epsilon$ exists since $x_p\to \widetilde p$ and 
$y_q\to \widetilde q$ as $r,u\to 1$.
In this case $$\frac{C}{||x_p - y_q||^{n-1}} \leq \Big(\frac{2}{d}\Big)^{n-1},$$ the volume $D$ of $\Omega$ is
finite, $|T|^{-1}\Vol(\sigma_{i\delta}(p,T))\leq 1/T_i,$
and  $|S|^{-1}\Vol(\sigma'_{j\varepsilon}(q,S)) \leq 1/S_j$ so the limit of $L$ is zero.

\medskip

Case 2. $r\in [0,1-\epsilon]$.

\begin{eqnarray*}L'&=&\int_{p\in\Omega}\int_{q\in\Omega}\int_{y_q\in\sigma'_{j\varepsilon}(q,S)}\int_{x_p\in\sigma_{i\delta}(p,T)}
\frac{1}{||x_p-y_q||^{n-1}}dy_qdr dt^idpdq\\
&=& \int\int\int\int_{t^i\in T^i}\int_{r=0}^{1-\epsilon}\frac{1}{||(1-r)\Phi(t^{i\delta},p) + r\widetilde p - y_q||^{n-1}}dy_qdu dt^idpdq
\end{eqnarray*}

Then  $\Phi(t^{i\delta},\cdot) = \Phi_{t^{i\delta}}$ is a volume-preserving diffeomorphism
of $\Omega$, so we can make the substitution $p'=\Phi(t^{i\delta},p)$ and get
$$ L'=\int\int\int\int_{r\in[0,1-\epsilon]}\int_{p\in\Omega}
\frac{1}{||(1-r)\Phi(t^{i\delta},p) + r\widetilde p - y_q||^{n-1}}dy_q dt^id\lambda(p)dqdr dt^i$$ 
$$= \int\int\int\int_{r\in[0,1-\epsilon]}\int_{p'\in\Omega}
\frac{1}{||(1-r)p' + r\widetilde p  - y_q||^{n-1}}dy_q dt^id\lambda(p')dqdr dt^i.$$
For each $r$ we let $p_r = (1-r)p' + r\widetilde p$.
Then $d\lambda(p_r) = (1-r)^nd\lambda(p')$ and $\Omega$ is replaced by
by $\Omega_r \subset \Omega$ (a contraction moving towards $\widetilde p$), so
$$L'= \int\int\int\int_{r\in[0,1-\epsilon]}\frac{1}{(1-r)^{n}}\int_{p_r\in\Omega_r}
\frac{1}{||p_r-y||^{n-1}}dr d\lambda(p_r)$$
$$\leq \int_{r\in[0,1-\epsilon]}\frac{1}{\epsilon^{n}}\int_{p_r\in\Omega}
\frac{1}{||p_r-y||^{n-1}}dr d\lambda(p_r)$$
$$\leq \int_{r\in[0,1-\epsilon]}\frac{\Gamma}{\epsilon^n}dr =\frac{\Gamma(1-\epsilon)}{\epsilon^n}
$$
by Lemma \ref{Gamma}, since $1-r\geq \epsilon$ and $\Omega_r\subset\Omega$. Now the volume of $\Omega$ is
finite, $\Vol(\sigma_{i\delta}(p,T))\leq |T|/T^i,$
and  $\Vol(\sigma'_{j\varepsilon}(q,S)) \leq |S|/S^j$,
so it follows that $\lim_{S,T\to\infty}L=0$.

\medskip

Case 3. $u\in [0,1-\epsilon]$. This case is exactly parallel to Case 2, with $p$ and $q$ interchanged, so it is omitted.
There is an overlap in the three cases, but all values of $(r,u)\in[0,1]\times[0,1]$ are covered.
\qed

\medskip

Then for almost all $(p,T)\in \Omega\times \T_k$ and $(q,S)\in
\Omega\times \T_\ell$, $\theta_\Phi(p,T)$ and
$\theta'_\Psi(q,S)$ are disjoint and the linking
number ${\rm
lk}(\theta_\Phi(p,T),\theta'_\Psi(q,S))$
is defined.

\bprop The limit 
\begin{eqnarray}\widetilde{\rm lk}(p,q)=\lim_{T_1,\dots,
T_k,S_1,\dots, S_\ell\rightarrow\infty} \frac{1}{T_1\dots
T_kS_1\dots S_\ell}{\rm
lk}(\theta_\Phi(p,T),\theta'_\Psi(q,S)) \label{tildelk}
\end{eqnarray}
exists as an integrable $L^1$-function on $\Omega\times\Omega$ and
does not depend on the choice of the points $\tilde p$ and $\tilde q$. \label{L1function}\eprop

\bpf.  
Calculating the linking number using the linking form (\ref{linkingform}),
it suffices to integrate over the sets $\Phi(p,T)$ 
and $\Psi(q,S)$, since by Proposition \ref{smallnbds} the limits of the integrals over the other three sets vanish, i.e.,

\begin{eqnarray}
\lim_{T,S\to\infty}\frac{1}{\lambda_k(T)\lambda_\ell(S)}\lk(\theta_\Phi(p,T),\theta'_\Psi(q,S)) =\nonumber\\
\lim_{T,S\to\infty}\frac{1}{\lambda_k(T)\lambda_\ell(S)}
\int_{\Phi(T,p)}\int_{\Psi(S,q))}L.
\label{fX}
\end{eqnarray}
As before, $X=X_1\wedge\dots\wedge X_k$ and $Y=Y_1\wedge\dots\wedge
Y_\ell$ are the exterior
products of the vector fields that generate the actions of $\Phi$
and $\Psi$, respectively. Define the function
$f:\Omega\times\Omega\rightarrow {\mathbb R}$ by
\begin{equation}f(p,q):=\frac{(-1)^k[(q-p)\times X(p)]\cdot Y(q)}{a_n||q-p||^n}\label{deffpq}.
\end{equation} 

For every $(p,q)$ we have
\begin{eqnarray}
|f(p,q)| \leq \frac{||X(p)||\ ||Y(q)||}{a_n||q-p||^{n-1}}
\leq \frac{K}{a_n||q-p||^{n-1}}\label{fpqleq}
\end{eqnarray}
where $K$ is an upper bound for $||X(p)||\hspace{0.2em} ||Y(q)||, p,q\in\Omega$.
Now,  by Lemma \ref{smallnbds}, $g(p,q)=1/||q-p||^{n-1}$ is an integrable function in $\Omega\times\Omega$, since
\begin{eqnarray*}
\int\int_{(p,q)\in\Omega\times\Omega}g(p,q)d\lambda(p)d\lambda(q) =
\int_{p\in\Omega}\left[\int_{q\in\Omega}\frac{1}{||q-p||^{n-1}}d\lambda(q)\right]d\lambda(p)\\
\leq \int_{p\in\Omega}\Gamma d\lambda(p) = \Gamma \Vol(\Omega).
\end{eqnarray*} Then by (\ref{fpqleq}) we get

\begin{eqnarray*}
\int\int_{(p,q)\in\Omega\times\Omega}|f(p,q)|d\lambda(p)d\lambda(q)  \leq \frac{\Gamma K\Vol(\Omega)}{a_n}
\end{eqnarray*} so $f\in L^1(\Omega\times\Omega)$.

To calculate $\int_{\Phi(T,p)}\int_{\Psi(S,q))}L$ we use the natural parametrizations $\bar p = x_p(t)=\Phi_t(p)=\Phi(t,p)$ and $\bar q = y_q(s)=\Psi_s(q)=\Psi(s,q)$ induced by the actions $\Phi$ and $\Psi$ on $\Phi(T,p)$ and $\Psi(S,q)$.
Then
$\frac{\pd x_p}{\pd t_i}(t) = X_i(\Phi_t(p)), i=1,\dots,k,$  
and $\frac{\pd y_q}{\pd s_j}(s) = Y_j(\Psi_s(q)), j=1,\dots,\ell$.
Let $\frac{\pd x_p}{\pd t} = \frac{\pd x_p}{\pd t_1}\wedge\dots\wedge   \frac{\pd x_p}{\pd t_k}(t)$ and
$\frac{\pd y_q}{\pd s}(s) = \frac{\pd y_q}{\pd s_1}\wedge\dots\wedge  \frac{\pd y_q}{\pd s_\ell}(s),$
so
$$\frac{\pd x_p}{\pd t}(t) = X_1\wedge\dots\wedge X_k(\Phi_t(p))  = X(\Phi_t(p))$$ and
$$\frac{\pd y_q}{\pd  s}(s) = Y_1\wedge\dots\wedge Y_\ell(\Psi_s(q) = Y(\Psi_s(q)).$$
Let $U(\bar p)$ and $U'(\bar q)$ denote the unit $k$- and $\ell$-vectors at $\bar p\in\Phi(T,p)$ and
$\bar q\in\Psi(S,q)$, respectively. Then by (\ref{linkingform})

\begin{eqnarray*}
(-1)^ka_n\int_{\Phi(T,p)}\int_{\Psi(S,q))}L \! =\! \int_{\bar p\in\Phi(T,p)}\int_{\bar q\in\Psi(S,q)}
\frac{[(\bar q-\bar p)\times U(\bar p)]\cdot U'(\bar q)}{||\bar q - \bar p||^n}d\eta(\bar p)d\eta(\bar q)
\end{eqnarray*}\begin{eqnarray*}
&=& \int_{t\in T}\Big(\int_{s\in S}  
\frac{[(y_q(s)-x_p(t))\times U(x_p(t))]\cdot U'(y_q(s))}{||y_q(s)-x_p(t)||^n}
||\frac{\pd y_q}{\pd  s}(s)||ds\Big)||\frac{\pd x_p}{\pd t}(t)||dt\\
&=&\!\!\int_{T}\int_{S}
\frac{[(\Psi_s(q)-\Phi_t(p))\times(||X(\Phi_t(p))||U(\Phi_t(p))]\cdot(||Y(\Psi_s(q))||U'(\Psi_s(q))}
{||\Psi_s(q)-\Phi_t(p)||^n}dsdt\\
&=&\int_{t\in T}\int_{s\in S}\frac{[(\Psi_s(q)-\Phi_t(p))\times X(\Phi_t(p))]\cdot Y(\Psi_s(q))}
{||\Psi_s(q)-\Phi_t(p)||^n}dsdt\\
&=&\int_{t\in T}\int_{s\in S}f(\Phi_t(p),\Psi_s(q))dsdt\\
&=&\int_{t\in T}\int_{s\in S}f(\Theta_{(t,s)}(p,q))ds_1\dots ds_\ell dt_1\dots dt_k\\
&=&
\int_0^{T_1}\dots\int_1^{T_k}\int_0^{S_1}\dots\int_0^{S_1}\dots\int_0^{S_\ell}
f(\Theta_{t_1,\dots,t_k,s_1,\dots,s_\ell}(p,q))ds_1\dots ds_\ell dt_1\dots dt_k,
\end{eqnarray*}
where $\Theta=\Phi\times\Psi$ is the product action of $\R^{k+\ell}$
in $\Omega\times\Omega$ defined by setting $\Theta_{(t,s)}(p,q) = (\Phi_t(p),\Psi_s(q))$.
Then the Ergodic Theorem, Theorem  \ref{ergthm}, applied to the action $\Theta$, shows that the limit
\begin{eqnarray*}
\lim_{T_1,\dots T_k,S_1,\dots S_\ell\to\infty} \frac{1}{T_1\dots T_kS_1\dots S_\ell}
\int_{\Phi(T,p)}\int_{\Psi(S,q))}L
\end{eqnarray*}
converges and defines a function $\widetilde{\lk}\in L^1(\Omega\times\Omega)$,
\begin{eqnarray*}
\widetilde{\lk}(p,q) &= \lim_{T_1,\dots T_k,S_1,\dots S_\ell\to\infty} \frac{1}{T_1\dots T_kS_1\dots S_\ell}\int_{\Phi(T,p)}\int_{\Psi(S,q))}L =
\end{eqnarray*}
\begin{eqnarray*}
\lim_{T,S\to\infty} \frac{1}{\lambda(T)\lambda(S)}
\int_0^{T_1}\!\!\!\dots\!\int_1^{T_k}\!\!\!\!\int_0^{S_1}\!\!\!\dots\int_0^{S_1}\dots\int_0^{S_\ell}\!\!
f(\Theta_{(t,s)}(p,q))dt_1\dots dt_kds_1\dots ds_\ell
\end{eqnarray*}
so we get
\begin{eqnarray*}
\int\int_{(p,q)\in\Omega\times\Omega}\widetilde{\lk}(p,q)dp\times dq = 
\int\int_{(p,q)\in\Omega\times\Omega}f(p,q)dp\times dq.
\end{eqnarray*}
Then (\ref{fX}) shows that this function satisfies (\ref{tildelk}). Clearly
it does not depend on the choices of $\widetilde p$ and $\widetilde q$.
\qed

\medskip

As a consequence of this Proposition, we can define the 
asymptotic linking invariant to be
$${\rm lk}(\Phi,\Psi)=\int_{p\in \Omega}\int_{q\in
\Omega}\widetilde{\rm lk}(p,q)d\eta(p)d\eta(q),$$ and then Theorem
\ref{mainthm} states that ${\rm lk}(\Phi,\Psi)=I(\Phi,\Psi)$.

\vskip 0.5cm
\noindent{\bf Proof of Theorem \ref{mainthm}.} 
With the volume forms $\omega, d\eta(p),$ and $d\eta(q)$ on $\Omega$, we have
\begin{eqnarray*}
{\rm lk}(\Phi,\Psi) &=&\int_{\Omega\times\Omega}\widetilde{\rm lk}(p,q)d\eta(p) d\eta(q)
= \int_{\Omega\times\Omega}f(p,q)d\eta(p) d\eta(q)\\
&=&
\frac{(-1)^k}{a_n}\int_{p\in\Omega}\int_{q\in\Omega}
\Big[\frac{q-p}{||q-p||^n}\times X(p)\Big]\cdot Y(q)d\eta(p)d\eta(q)
\quad\quad\quad {\textnormal  by\ (\ref{deffpq})}\\
&=& \int_{q\in\Omega}\Big[\frac{(-1)^k}{a_n}\int_{p\in\Omega}
\frac{q-p}{||q-p||^n}\times X(p)d\eta(p)\Big]\cdot Y(q)d\eta(q)
\end{eqnarray*}
by Fubini's Theorem, and then, by the Biot-Savart formula  (\ref{eqBS}), the definition of $j$,
 (\ref{alphaY}), and Corollary \ref{c-djBS}, this is equal to
\begin{eqnarray*}
\int_\Omega (BS(X)\cdot Y)\omega
= \int_\Omega jBS(X)(Y)\omega 
\end{eqnarray*}
 $$=  \int_\Omega jBS(X)\wedge i_Y\omega = \int_\Omega jBS(X)\wedge d\beta = I(\Phi,\Psi).$$  
\qed


\section{A lower bound for the energy of an action}\label{energy} 

We remark that in the case when $\Phi=\Psi$ and $n=2k+1$, the
invariant ${\rm lk}(\Phi,\Phi)=I(\Phi,\Phi)$ is a lower bound for
the energy of the generating $k$-vector $X$. 

\bdf Let  $\Phi$ be a conservative $k$-action on $\Omega$ and let $X$ be the $k$-vector field that generates $\Phi$.
The {\em energy} of the $k$-action $\Phi$ is defined to be the value of the integral
$$E(\Phi) = ||X||^2 = \int_{p\in \Omega}X(p)\cdot X(p)d\lambda(p) = 
\int_{p\in \Omega}||X(p)||^2d\lambda(p).$$
\edf

Note that we can decrease  the energy of $\Phi$ by conjugating $\Phi$ by volume-preserving diffeomorphisms.
Can we make it arbitrarily close to zero? The following result gives a negative answer to this question.
\bthm There exists a constant $C>0$ depending only on $\Omega$ such that
$$C^{-1}|{\rm lk}(\Phi,\Phi)| \leq E(\Phi).$$
\ethm
\bpf. By Corollary 2, (7),  and the definition of $j$,
\begin{eqnarray*}{\rm lk}(\Phi,\Phi) &=& \int_\Omega jBS(X)\wedge d\alpha = 
\int_\Omega jBS(X)\wedge i_Xd\lambda\\
&=& \int_\Omega jBS(X)(X) d\alpha = 
\int_\Omega BS(X)\cdot X d\lambda.
\end{eqnarray*}
By the Cauchy-Schwarz inequality
\begin{eqnarray}|{\rm lk}(\Phi,\Phi)| = |<\!BS(X),X\!>| \leq ||BS(X)||\ ||X||.\label{CauchySch} \end{eqnarray}
Furthermore
$$BS(X)(p) = \int_{q\in\Omega}\frac{(p-q)\times X(q)}{||p-q||^{2k+1}}d\lambda(q)$$
so
\begin{eqnarray*}||BS(X)(p)|| &\leq& \int_{q\in\Omega}\frac{||(p-q)\times X(q)||}{||p-q||^{2k+1}}d\lambda(q)\nonumber\\
&\leq& \int_{q\in\Omega}\frac{(||X(q)||}{||p-q||^{2k}}d\lambda(q))\nonumber\\
&=& \int_{q\in\Omega}\Big[\frac{||X(q)||}{||p-q||^{k}}\Big]\Big[\frac{1}{||p-q||^k}\Big]d\lambda(q))\nonumber\\
&\leq& \Big[\int_{q\in\Omega}\frac{||X(q)||^2}{||p-q||^{2k}}d\lambda(q)\Big]^{1/2}\Big[\int_{q\in\Omega}\frac{1}{||p-q||^{2k}}\lambda(q)\Big]^{1/2}\end{eqnarray*}
by the Holder inequality. Then by Lemma \ref{Gamma} with $n=2k+1$
\begin{eqnarray}||BS(X)(p)|| \leq \Gamma^{1/2}\int_{q\in\Omega}\Big[\frac{||X(q)||^2}{||p-q||^{2k}}d\lambda(q)\Big]^{1/2}.\label{ineq2}\end{eqnarray}
 Therefore

\begin{eqnarray*}
||BS(X)||^2 &=& \int_{p\in\Omega}BS(X)(p)\cdot BS(X)(p)d\lambda(p)\nonumber\\
&=& \int_{p\in\Omega}||BS(X)(p)||^2d\lambda(p)\nonumber\\
&\leq& \Gamma\int_{p\in\Omega}\Big[\int_{q\in\Omega}\frac{||X(q)||^2}{||p-q||^2}d\lambda(q)\Big]d\lambda(p)
\quad\quad{\textnormal \rm by\ \  (\ref{ineq2})}\\
&=& \Gamma\int_{q\in\Omega}||X(q)||^2\Big[\int_{p\in\Omega}\frac{1}{||p-q||^2}d\lambda(p)\Big]d\lambda(q)\ \ \textnormal{by Fubini's\ Theorem}\\
&\leq& \Gamma^2\int_{q\in\Omega}||X(q)||^2d\lambda(q)
\end{eqnarray*}
by Lemma \ref{Gamma}. Thus 
$$ ||BS(X)|| \leq (\Gamma^2)^{1/2}\Big)\int_{q\in\Omega}||X(q)||^2d\lambda(q)\Big)^{1/2} = \Gamma||X||.$$

Substituting this inequality in (\ref{CauchySch}) we get
\begin{eqnarray}
|{\rm lk}(\Phi,\Phi)| \leq ||BS(X)||\ ||X|| \leq \Gamma||X||^2| = \Gamma E(\Phi). \label{energybd}
\end{eqnarray} 

We can decrease the energy of $\Phi$ by volume-preserving diffeomorphisms, but these diffeomorphisms
do not change the value of the asymptotic linking number ${\rm lk}(\Phi,\Phi)$, so by (\ref{energybd}) $\Gamma^{-1}{\rm lk}(\Phi,\Phi)$ is the desired lower bound for the energy of $\Phi$.\qed

\section{Examples}\label{examples} 

\bex For every pair of integers $k,\ell \geq 1$, $k+\ell +1 = n$, and every $t\in\R$, there are conservative actions 
$\Phi$ of $\R^k$ and $\Psi$ of $\R^\ell$ on the unit closed ball $D^n\subset \R^n$
such that ${\rm lk}(\Phi,\Psi)=I(\Phi,\Psi)=t$. \label{ex1}\eex

The construction uses several lemmas.

\blm Given disjoint smooth embeddings of closed oriented manifolds $M, N$, of dimensions $k$ and $\ell$ 
in ${\mathbb R}^n$, there exist disjoint smooth embeddings $M\times S^1, N\subset{\mathbb R}^{n+1}$
such that ${\rm lk}(M\times S^1,N) =  {\rm lk}(M,N)$. The same holds if $N$ is an affine $\ell$-space
disjoint from $M$. \label{multS1}\elm

\bpf. Given $M$ and $N$, by a translation we may assume that their images lie in the positive half space
 $x_1>0\subset {\mathbb R}^n\subset{\mathbb R}^{n+1}$. Let  $P$ be the $(n-1)$-plane  in ${\mathbb R}^n$
 perpendicular to the $x_1$-axis, and  
 rotate $M$ around $P$ to get $M\times S^1\subset {\mathbb R}^{n+1}$.
 Clearly $M\times S^1$  is disjoint from $N$. If we let $\Sigma\subset  {\mathbb R}^n$ be a compact singular $(k+1)$-manifold
tranverse to $N$ such that $\partial\Sigma = M$, then ${\rm lk}(M,N) = {\rm Int}(\Sigma,N)$.
By rotating $\Sigma$ around $P$ we obtain $\Sigma\times S^1$, whose boundary is $M\times S^1$. Then
$ {\rm Int}(\Sigma\times S^1,N) = {\rm Int}(\Sigma,N)$, and therefore the linking number is the same.

In case $N$ is an affine $\ell$-plane a similar argument works, taking $P$ to be an affine plane parallel to $N$.
 \qed

\blm There exist disjoint embeddings of $T^k\times D^{\ell+1}$ and $T^\ell\times D^{k+1}$ in $D^n$, where $T^k$ and $T^\ell$ are tori of dimensions $k$ and $\ell$, such that
$\rm{lk}(T^k\times  0, T^\ell\times 0)=1$. \label{Tk} \elm

\bpf. Begin with disjoint smooth embeddings of two circles $M$ and $N$ in ${\mathbb R}^3$ such that
${\rm lk}(M,N) = 1$. Applying Lemma \ref{multS1} repeatedly, switching the roles of $M$ and $N$, gives
disjoint embeddings of $T^k$ and $T^\ell$ in ${\mathbb R}^n$ with intersection number $1$. Since the normal bundles are trivial
we can extend the embeddings to disjoint embeddings of  $T^k\times D^{\ell+1}$ and $T^\ell\times {k+1}$.
then a homothety will move these sets into $D^n$. \qed

\blm Let $T^r$ act on $T^r\times D^s$ by the product action on the first factor
and identity on the second factor.
For any smooth volume form $\omega$ on $T^r\times D^s$ there is a
smooth isotopy 
$h_t$ of $T^r\times D^s$ taking each factor $T^r\times \{y\}$ to itself such that 
$h_0={\rm id}$ and $h_1^*(\omega)$ is $T^r$-invariant.
\label{isotopy}\elm

\bpf. Here we need a slightly modified form of Moser's Theorem \cite{M} acting on each orbit. We use the standard coordinates  $(x,y)=(x_1,\dots,x_k,y_1,\dots, y_s)$  and the standard Euclidean volume form 
$\omega^*=dx\wedge dy$ on $T^r\times D^s$ to simplify the notation.
Let $f_0: T^r\times D^s\to\R$ be the (unique)  non-vanishing smooth function such that $\omega=f_0\omega^*$
and define $f_1:D^s\to\R$ by setting $f_1(y) = \int_{T^k\times \{y\}}f_0(x,y)dx$.
Note that $\alpha^*_y = f_0(y)dx$ is $T^k$-invariant.
Now the volume forms $\alpha_y = f_1(x,y)dx$ and $\alpha^*_y$ have the same integral
$\int_{T^k}\alpha_y dx = \int_{T^k}\alpha^*_y dx$, so there exists a smooth function
$f: T^r\times D^s\to\R$ such that $\alpha=f\alpha^*$ and
we can apply Moser's proof  \cite{M} on each factor $T^k\times \{y\}$.
Following Moser, we may suppose that there is a positive $\epsilon$ such that $|f(y)-1|<\epsilon$ for every $y$
by expressing any positive function $f$ as a sum of functions close to $1$.
We use the same cover of $T^k$ by open cubes $U_0,U_1,\dots,U_m$ and the same functions
$\eta_i (k=1,\dots, m)$, independent of $y$. Then it is straightorward to check that Moser's isotopies
of each $T^k\times \{y\}$ fit together to give a smooth isotopy of $T^r\times D^s$ transforming
each $\alpha_y$ into $\alpha_y^*$. This isotopy also transforms $\omega$ into $f_1\omega_*$, which is 
invariant under the action of $T^k$. \qed

\bigskip

\noindent{\bf Construction of the Example.} Take
$W=T^k\times D^{\ell+1} \sqcup T^\ell\times D^{k+1}$ embedded
in $D^n$ by  Lemma \ref{Tk}, where $k+\ell+1=n$.
The compact Lie groups $T^k$ and $T^\ell$ act on $W$, 
$T^k$ acting on $T^k\times D^{\ell+1}$ by multiplication on the first
factor and trivially on $T^\ell\times D^{k+1}$, and analogously for the
action of $T^\ell$. 

By Lemma \ref{isotopy},
we may conjugate the action of $T^k$ on $T^k\times D^{\ell+1}$ by a diffeomorphism
isotopic to the identity so that it preserves
the Euclidean volume form, and similarly for the action of $T^{\ell}$. 
Lift the actions of $T^k$
and $T^\ell$ to volume preserving actions  $\phi:\R^k\times
(T^k\times D^{\ell+1})\to T^k\times D^{\ell+1}$ and 
$\psi:\R^\ell\times (T^\ell\times D^{k+1})\to T^\ell\times D^{k+1}$.

Let $W_\epsilon =  T^k\times D_0^{\ell+1} \cup T^\ell\times D_0^{k+1}$ be a smaller
invariant neighborhood of $T^k \cup T^\ell$  and 
let $\lambda:W\to [0,1]$ be constant on the orbits with the values $1$ on $W_\epsilon$ and $0$
on $D^n\setminus W$. Then let $\Phi(t,z)=\phi(\lambda(z)t,z)$ for $z$ in the $\epsilon$-neighborhood of
$T^k\times D^{\ell+1}$ and identity elsewhere, while $\Psi(t,z)=\psi(\lambda(z)t,z)$
on the $\epsilon$-neighborhood of $T^\ell\times D^{k+1}$ and identity elsewhere.
Thus $\Phi$ and $\Psi$ are commuting conservative actions of $\R^k$ and $\R^\ell$ on $D^n$.
The linking number of the orbits $T^k\times \{y\}$ and $T^\ell\times \{z\}$ are 
$\rm{lk}(T^k\times  y, T^\ell\times z)=1$ for $y\in D_0^{\ell+1}$ and $z\in D_0^{k+1}$.

Now it is easy to check that the linking number $\rm{lk}(\Phi,\Psi) > 0$
since for points $p\in D_0^{\ell+1}$ and $q\in D_0^{k+1}$
and for $T=[0,2r\pi]^k$ and $S=[0,2s\pi]^\ell$,
$$\rm{lk}(\theta_\Phi(p,T),\theta_\Psi(q,S)= 
r^ks^\ell$$ 
since for these rectangles $T$ and $S$ the cones $\sigma(p,T)$ and $\sigma'(q,S)$ are empty.
When we normalize by dividing by $(2r\pi)^k\cdot (2s\pi)^\ell$
we get the constant $(2\pi)^{-(k+\ell)}$, which is therefore the value of the limit 
for orbits in $W_\epsilon$ as $r,s\to\infty$.
Other points $p, q$ contribute positively, so we get
$\rm{lk}(\Phi,\Psi) > 0$. To get a negative value it suffices to change one
of the orientations. Finally by multiplying $t\in \R^k$ by $s$
we multiply the asymptotic linking number by $s^k$ and thus we can obtain all real numbers
as values of $\rm{lk}(\Phi,\Psi)$. \qed

\bex Given a closed connected oriented submanifold $N^\ell$ embedded in $D^n$ 
and a real number $t$, by a similar construction we can find a conservative action 
$\Phi$ of $\R^k$ on $D^n$, $k=n-\ell-1$, such that $\rm{lk}(\Phi,N)=t$. 
\eex

Here the construction is similar to the previous example. By applying Lemma \ref{multS1}
repeatedly we can obtain $T^k\subset \R^n\setminus P$, where $P$ is an affine $\ell$-plane,
such that the linking number is $\rm{lk}(T^k,P) = 1$. Now locally the smooth embedding
of $N$ in $D^n$ is diffeomorphic to the embedding of $P$ in $\R^n$, so we can find a
small torus $T^k\subset D^k\setminus N$ such that $\rm{lk}(T^k,N) = 1$. The
rest of the construction proceeds as in Example \ref{ex1}.

\end{document}